\documentstyle[12pt]{article}
\input amssym.def
\input amssym

\newtheorem{Theo}{Theorem}
\newtheorem{Lemm}{Lemma}[section]
\newtheorem{Prop}{Proposition}
\newtheorem{Cor}{Corollary}
\newcommand{\pf}{{\bf Proof}\ }

\newcommand{\Q}{{\bf Q}}
\newcommand{\Z}{{\bf Z}}
\newcommand{\N}{{\bf N}}
\newcommand{\C}{{\bf C}}
\newcommand{\D}{{\bf D}}
\newcommand{\E}{{\cal E}}
\newcommand{\EZ}{{\cal S}}
\newcommand{\LL}{{\cal L}}
\newcommand{\spn}{{\rm span}}
\newcommand{\ssp}{{\rm span}_Q}
\newcommand{\EC}{{\cal SC}}
\newcommand{\SED}{sub{\cal S}}
\newcommand{\SEO}{sub{\cal E}}
\newcommand{\nt}{{\bf Notation}\ }

\newcommand{\card}{{\rm card\ }}
\newcommand{\trd}{{\rm tr.d. }}
\newcommand{\ldim}{{\rm tr.d. }}
\newcommand{\adim}{{\rm adim\ }}
\newcommand{\Ker}{{\rm Ker}}
\newcommand{\mdim}{{\rm dim_Q }}
\newcommand{\ex}{{\rm ex}}
\newcommand{\ecl}{{\rm ecl}}
\newcommand{\acl}{{\rm acl}}
\newcommand{\cl}{{\rm cl}}

\newcommand{\ssn}{\section}
\newcommand{\rem}{{\bf Remark }}
\newcommand{\df}{{\bf Definition}\ }
\newcommand{\bl}{\begin{Lemm}}
\newcommand{\el}{\end{Lemm}}
\newcommand{\bt}{\begin{Theo}}
\newcommand{\et}{\end{Theo}}

\newcommand{\bp}{\begin{Prop}}
\newcommand{\ep}{\end{Prop}}

\newcommand{\bc}{\begin{Cor}}
\newcommand{\ec}{\end{Cor}}
\newcommand{\lb}{\label}
\newcommand{\dd}{\partial}

\newcommand{\ra}{\rangle}
\newcommand{\la}{\langle}
\newcommand{\subs}{\subseteq}
\newcommand{\qed}{$\Box$ \\ \\}
\newcommand{\smin}{\setminus}
\newcommand{\sups}{\supseteq}
\begin{document}
\title{Fields with pseudo-exponentiation}
\author{B. Zilber\\ zilber@maths.ox.ac.uk }  \date{13 October 2000}
\maketitle

\ssn{Introduction}

This research is motivated by the study of model-theoretical
properties of classical 'analytic' structures, i.e. ones having
natural analytic representation (see also [Z]). For example, the
structure of complex
numbers as a field with exponentiation 
$$\C_{\exp} =(\C,+,\cdot,\exp).$$
One of the questions we can ask is whether $\C_{\exp}$ is
{\bf quasi-minimal,} i.e. any definable
subset of $\C_{\exp}$  is either countable or 
of power continuum. Another question is about homogeneity of the
structure; we do not know any its automorphism except the identity and
the complex conjugation.
In general we would like to understand the nature of analytic
dimension in a context close to model-theoretic stability theory.

A slightly weaker analytic structure $\C^{(2)}_{\exp}$ is a
two-sorted structure with both sorts $\C(1)$ and $\C(2)$ being copies
of complex numbers, on both sorts the field structure is given and
there is a mapping  $\exp:\C(1)\to \C(2)$ in the language.

The model theory of the both structures, as well as of many others of
this 
kind, seems very hard to study directly. 

We study here model theory of abstract analogues of $\C^{(2)}_{\exp}.$
We start by considering
 the class $\E_p$ of two-sorted structures of the form $(D,\ex,R),$ \ 
$D,$ the domain of a mapping $\ex$ is a field of  characteristic zero,
$R,$ a field of  characteristic $p$      and $\ex $ is a homomorphism of the additive group of $D$ onto the multiplicative group $R^{\times}$ of $R.$
Following the ideology of Hrushovski's construction of non-classical
structures with nice dimension notion (see [H]), we consider a notion of a {\em predimension} $\delta(X)$ for finite subsets
$X\subs D$ defined as
$$\delta(X)=\ldim(X)+\trd (\ex (X))- \mdim(X),$$ here $\mdim(X)$ is the
dimension of the $\Q$-linear space generated by $X.$\\
We then define the subclass $\EZ_p$ of $\E_p,$ the {\em weak Schanuel class}
  as the class of structures where $\delta(X)\ge 0$ for all $X$ and the
kernel of $\ex$ is {\em standard} (a cyclic additive subgroup in case  $p=0$).
Recall that the Schanuel  conjecture [L] states
that, given $\Q$-linearly independent complex numbers $x_1,\dots,x_n,$ the
inequality  $$\trd(x_1,\dots,x_n, \exp(x_1),\dots,\exp(x_n))-n\ge 0$$ 
holds.  We prove that the weak Schanuel  class $\EZ_p$ is non-empty.

  We study the predimension $\delta$ and the
relation $A\le B$ following the Hrushovski
pattern.  We define a notion
of {\em exponentially-algebraically closed structures}, which are just
existentially closed structures  with respect to the relation $A\le B$ in class
$\EZ_p.$ 
Again, by standard model theory as used by Hrushovski, in the class $\EC_p$ of
exponentially-algebraically closed  structures $\delta$ gives rise 
 to the {\em dimension notion} $\partial.$ 

On the other hand class $\EC_p$ is shown to be axiomatizable inside
$\EZ_p.$ More precisely, it is axiomatized by an axiom stating that
any
  {\em normal}
system of exponential-algebraic equations has a solution in the structure.  
The definition of a  normal
system of exponential-algebraic equations is explicit and purely
algebraic. 
 It follows
  
 {\em A necessary and sufficient condition for the
structure of complex numbers with $\exp$ to be in $\EC_0$ is to satisfy
the weak Schanuel conjecture and to have solutions to any  normal
system of exponential-algebraic equations.}  

There are grounds to believe that complex numbers do satisfy the axiom.\\

Next we proceed by developing analogues of atomic and prime structures
 in class $\EC_p.$ A notion
of quasi-atomic types and  structures as well as quasi-prime structures
over a given subset are introduced in style of [Sh]. We restrict the study to the case
$p=0.$ 
Then the usual nice properties for the
notions hold. Given a subset $C\le \D$ in a structure in
$\D\in\EC_0$ we construct a quasi-prime
structure $E(C)\le \D.$ If the kernel of $\ex,$ the additive subgroup
$K=\{ x\in D: \ex(x)=1\},$ is  compact in the profinite topology then $E(C)$  is determined uniquely up to isomorphism over $C.$ Now
take $C$ a  $\partial$-independent set of power $\lambda.$ We call
$E(C)$ in this case {\em the canonical weak Schanuel structures
of power  $\lambda$} (canonical for short). This notion may be useful
because  $\C^{(2)}_{\exp}$ looks very much like a
canonical structure of power $2^{\aleph_0}$ with the {\em standard} kernel
 isomorphic to $\Z.$

 This is certainly hard to prove, since
the proof would include a proof of the weak Schanuel conjecture.  Still it
would be interesting to study the properties of  canonical
structures compairing these to the known properties of $\C.$   It follows
from what was stated above that  

{\em For
$\lambda>\aleph_0,$ a canonical structure with the standard kernel is  
quasi-minimal, i.e. any definable
subset of $E_{\lambda}$ is either  of power not greater than $\aleph_0$ or
of power $\lambda.$}\\  

In a canonical structure with standard kernel $E_{\lambda}$ ($\lambda>\aleph_0$) we study the
notion of dimension
for definable subsets $S\subs E^n,$ denoted $\adim(S),$  which we call
{\bf the pseudo-analytic dimension}. We define it analogously to the definition
of Morley Rank: $$\adim S=\max_{ \bar x\in S} \partial \bar x.$$ The
pseudo-analytic dimension  is a good analogue of the analytic dimension in
$\C:$  {\em

(i) $\adim S\ge 0$ for any non-empty $S;$

(ii) for non-empty $S$ \ $\adim S=0\mbox{ iff }\card (S)\le\aleph_0;$

(iii)  for non-empty $S$ \ $\adim S\ge m>0$ iff $S$ can be projected to
$D^m$ so that the complement to the projection in $D^m$ is of pseudo-analytic
dimension less than $m$;

(iv) For any irreducible algebraic varieties $V\subs D^n$ and $W\subs R^n$

\hspace{1cm}  $V\cap \ln W\neq\emptyset$ implies
$\adim(V\cap \ln W)\ge \dim V+\dim W-n.$}\\ \\

{\bf Acknowledgements}. The research had been conceived in February
1997 while the author was visiting The Fields Institute in Toronto. The study
continued in Kemerovo University, Russia, where I had
a permanent position then. Essential progress had been achieved while I
was in MSRI, Berkeley, participating in the Model Theory of Fields
program.  
The final version of the paper has been written in University of
Oxford. I would like to express my gratitude to all the
institutions and people there who supported and helped the research.

\ssn{Definitions and notation}

For technical reasons it is more suitable to work with the equivalent class
of one sorted structures on $D$ in the language $\LL_p$ which is the version  of extended language of fields. On $D$ the following operations and relations are defined:  the addition $' + '$ and multiplication by all rational scalars;
$n$-ary predicates $V(x_1,\dots,x_n)$ for each  algebraic subvariety
$V\subs D^n$ defined and irreducible over $0$
 an equivalence relation $E,$ this allows to consider $R=D/E,$ we think of $R$ as of a field;   $n$-ary predicates $EW(y_1,\dots,y_n)$ for each algebraic subvariety $W\subs R^n$ definable and irreducible over $0.$\\

The interpretation  assumes     $E(x,y)\equiv [\ex(x)=\ex(y)],$

   $V(x_1,\dots,x_n)\equiv [\la x_1,\dots,x_n\ra \in V],$

   $EW(x_1,\dots,x_n)\equiv [\la \ex(x_1),\dots,\ex(x_n)\ra \in W].$\\ \\
\df
$\E_p$ is the class of structures $\D$ in language $\LL_p$ such that
 $D$ is a field of characteristic zero,
 $E$ an equivalence relation on $D$ which is
congruent with respect to the relations $EW(x_1,\dots,x_n),$
$D/E$  when considered with relations $EW(x_1,\dots,x_n)$ is embeddable
in the group $R^{\times}$ where $R$ is a field
of characteristic $p$    so that the canonical mapping
$$\ex: D\to R^{\times} $$ is a homomorphism of the additive group of
$D$ into the multiplicative group $R^{\times}.$\\ \\
\nt
For the sake of simplicity we fix $p$ below and write $\E$ instead of
$\E_p.$\\

For finite $X, X'\subs D,$ \ $Y,Y'\subs R$

$\ldim(X) $ the
 the transcendence degree of $X$ over $0;$

$\mdim(X) $ the  dimension of the vector space $\spn_{\Q}(X)$ generated by
$X$ over $\Q;$

$\trd(Y)$ the transcendence degree of $Y$ over $0;$

$\delta(X)$ {\bf the predimension} of finite $X\subs D:$

$\delta(X)=\ldim(X)+\trd (\ex (X))- \mdim(X);$

$\delta(X/X')=\delta(XX') - \delta(X'),$ \ and \ $XX'=X\cup X';$\\

for infinite $A\subs D$ \ $\delta(X/A)\ge k,$ $k\in \Z,$ means that for
any finite $Y\subs A$ there is finite $Y\subs Y'\subs A$ such that
 $\delta(X/Y')\ge k,$ \ and $\delta(X/A)=k$ 
means $\delta(X/A)\ge k$ and not $\delta(X/A)=k+1.$  
 \\

$\ldim(X/X')=\ldim(XX')-\ldim(X');$

$\trd(Y/Y')=\trd(YY')-\trd(Y');$

$\mdim(X/X')=\mdim(XX')-\mdim(X');$

$\ker$ is the name of a unary predicate of type $EW$
$x\in \ker \equiv \ex(x)=1.$\\ \\
Denote $\EZ$ the subclass of $\E$ consisting of all $D$
satisfying the condition:
   $$\delta(X)\ge 0\mbox{ for all finite }X\subs D.$$
 Denote $\SEO$ the class of structures $A$ in language $\LL_p$
satisfying :

 $A$ is an additive divisible subgroup of the additive group of a field
$D$ of characteristic zero
and $A/E$ can be embedded into $R^{\times}$ for a field $R$ of characteristic $p$ and the
canonical mapping
$\ex_A:A\to R^{\times}$ is a homomorphism.\\ \\
 Denote $\SED$ the subclass of $\SEO$ consisting of $A$
which satisfy $\delta(X)\ge 0$ for any finite $X\subs A.$\\ \\
For $W$  an algebraic variety, $\bar b=\la b_1,\dots,b_l\ra$
denote
$$W(\bar b)=\{ \la x_{l+1},\dots,x_{n+l}\ra :   \la  b_1,\dots, b_l,  x_{l+1},\dots,x_{n+l}\ra\in W \}.$$ \bl \lb{q} If $X=\{ x_1,\dots x_{n+l}\}\subs D,$ \
$X'=\{ x_{1},\dots x_{l}\},$
$\bar x=\la x_1,\dots x_{n+l}\ra,$ \   $\bar x'=\la x_{1},\dots x_{l}\ra$  \  then:  $\ldim(X)=\dim V,$ where $V\subs D^{n+l}$ is the minimal algebraic variety
over $0$   containing $\bar x;$

$\mdim(X)=\dim  L$ where $\bar L$ is the minimal linear subspace of $D^{n+l}$
containing $\bar x$ and given by homogeneous linear equations over $\Q;$

$\trd(\ex X)$ is the dimension of the minimal subvariety of $R^{n+l}$
over $0$ containing $\ex(\bar x);$

$\mdim(X/X')= \dim L(0^l),$ where $0^l$ is a string of $l$ zeroes.\el
\pf Immediate from definitions.\qed
\bl \lb{d}

(i) $\delta(X/X')=\ldim(X/X')+\trd(\ex X/\ex X')-\mdim(X/X');$

(ii) For any $A\subs D$ and finite $X\subs D$ there is a finite $Y\subs A,$
such that, if $Y\subs Y'\subs A$, then $\delta(X/Y')=\delta(X/A).$\el
\pf (i) is immediate from definitions.\\
(ii) follows from (i) if we choose $Y\subs A,$ such that
$\trd(X/Y)=\trd(X/A),$ $\trd(\ex(X)/\ex(Y))=\trd(\ex(X)/\ex(A)),$
$\mdim(X/Y)=\mdim(X/A).$ The choice is possible, since
$\trd(X/Y),$ $\trd(\ex(X)/\ex(Y))$ and $\mdim(X/Y)$ are non-increasing
functions of $Y.$ \qed

{\bf Remark} The condition $\mdim(X/A)=\mdim(X/Y)$ for $Y\subs A$ is
satisfied iff $\ssp(X)\cap A\subs Y.$ Correspondingly with the
transcendense degree.\\ \\

\nt For $A, B\in \SEO$ denote $A\le B$ the fact that
$A\subs B$ as structures and $\delta(X/A)\ge 0$
for all finite $X\subs B.$

\bl \lb{delta} For any structure $A$ of the class $\SEO$ and finite
$X,Y,Z\subs A:$

(i) If $\ssp(X')=\ssp(X)$ then
$\delta(X')=\delta (X).$

(ii) If $\ssp(X'Y)=\ssp(XY)$ then
$\delta(X/Y)=\delta (X'/Y).$

(ii) If $\ssp(Y)=\ssp(Y')$ then
$\delta(X/Y)=\delta (X/Y').$

(iv)
$\delta(XY/Z)=\delta (X/YZ)+ \delta(Y/Z).$
 \el
\pf Immediate from the definitions.$\Box$

\bl \lb{tr}           For $A,B,C\in \SEO$

(i)  if $A\le B$ and $B\le C,$ then $A\le C;$

(ii) if $A\le B,$ $Y\subs B,$ $\delta(Y/A)=0,$ then $AY\le B.$
\el
\pf (i) Let $ X\subs C$ finite and $Z\subs A$ large enough 
finite so that  $\delta(X/Z)= \delta(X/A) .$ We need to prove that $\delta(X/Z)\ge 0.$ Choose $Y\subs B$
finite so that $\ssp(YZ)=A\cap \ssp(XZ).$ Then
$\mdim(X/YZ)=\mdim(X/B).$ From the definition of $\delta$ it follows 
that $\delta(X/YZ)\ge \delta(X/B)\ge 0.$ Hence 
$\delta(XY/Z)=\delta (X/YZ)+ \delta(Y/Z)\ge 0.$ At last notice that
$\delta(X/Z)=\delta(XY/Z)$ by Lemma~\ref{delta} (ii).\\

(ii) For any $X\subs B$ we want to show $\delta(X/AY)\ge 0.$ It is given by
$$\delta(X/AY)=\delta(XY/A)-\delta(Y/A)=\delta(XY/A)\ge 0,$$ since $A\le B.$
\qed
\nt
 Let $A\in \SEO$ and
$X\subs A$ finite. Denote
   $$\partial_A(X)=\min\{\delta(X'): X\subs X'\subs A\}.$$

\bl  \lb{1} Let $A\in \SEO$ and $X\subs A$ finite. Choose $X'\subs A$ finite
 such that
   $$\delta(X')=\partial_A(X).$$
Then $X'\le A.$\el
\pf Immediately from the definition.\qed

\bl  \lb{X} Let $A, B\in \SEO,$  $A\le B$ and $X$ a finite subset of $A.$
Then $$\partial_A(X)=\partial_B(X).$$\el
\pf Let a finite $Y\subs B$ be such that $$\delta(X Y)=\partial_B(X).$$

Let $Y_0$ be a $\Q$-linear basis of $Y$ over $A$ and $X_1\subs A$ a finite
superset of  $X$ such that $$\spn_Q(X_1)=\spn_Q(X Y)\cap A.$$
Then
$\mdim(Y_0/A)= \mdim(Y_0/X_1),$ \
$\ldim(Y_0/A)\le \ldim(Y_0/X_1),$  \ and \
$\trd(\ex Y_0/\ex A)\le \trd(\ex Y_0/\ex X_1).$ It follows
$$\delta (Y_0/X_1)\ge \delta(Y_0/A)\ge 0.$$

Also $$\spn_Q(X Y) = \spn_Q(X_1 Y_0).$$
Hence $$\delta(X Y)=\delta(X_1 Y_0)=\delta(X_1)+\delta(Y_0/X_1).$$
By the above proved  $\delta(X Y)\ge\delta(X_1)$ and by the definition
$\delta(X_1)\ge \partial_A(X).$ Thus
$\partial_B(X)\ge  \partial_A(X),$ and the converse holds by the definition.\qed

\bl \lb{ext} Suppose $A\in \SED,$ $A'\in \SEO,$ 
$A'=\spn_Q(A X),$  and 
$\delta(X'/A)\ge 0$ for all finite $X'\subs \spn_QX.$
Then $A'\in \SED$ and $A\le A'.$ \el
\pf We may assume that
$X$ is  $\Q$-linearly independent over $A.$
Let $Z\subs A',$ \ $Z=\{ z_1,\dots z_n\},$ and $z_i=x_i+y_i$ for some
$x_i\in \spn_Q(X),$  \ $y_i\in A.$
Let $\{ x_1,\dots x_k\}$ be a $\Q$-linear basis of $\{ x_1,\dots x_n\}.$
Then, using Lemma~\ref{delta},

 $\delta(Z)=\delta(x_1+y_1,\dots x_k+y_k,y'_{k+1},\dots,y'_n)$ \
for $y'_{k+1},\dots,y'_n$  appropriate $\Q$-linear combinations of
$y_1,\dots y_n.$\\ Rewrite
 $$\delta(Z)=
\delta(\{ x_1+y_1,\dots x_k+y_k\} /\{ y'_{k+1},\dots,y'_n\})+
\delta(y'_{k+1},\dots,y'_n).$$

By the assumtions $\delta(y'_{k+1},\dots,y'_n)\ge 0.$ On the other hand

$\delta(\{ x_1+y_1,\dots x_k+y_k\} /\{ y'_{k+1},\dots,y'_n\})\ge
\delta(\{ x_1,\dots x_k\}/A)\ge 0$\\
since

$\ldim(\{ x_1+y_1,\dots x_k+y_k\}/\{ y'_{k+1},\dots,y'_n\})\ge
\ldim(\{ x_1+y_1,\dots x_k+y_k\}/A)\ge
\ldim(\{ x_1,\dots x_k\}/A),$

$\trd(\ex\{ x_1+y_1,\dots x_k+y_k\}/\ex \{ y'_{k+1},\dots,y'_n\})\ge$\\
$\trd(\ex \{ x_1+y_1,\dots x_k+y_k\} /\ex A)\ge
\trd (\ex\{ x_1,\dots x_k\} /\ex A)$\\
and

$\mdim(\{ x_1+y_1,\dots x_k+y_k\}/\{ y'_{k+1},\dots,y'_n\})=k=
\mdim(\{ x_1,\dots x_k\}/A).$
Thus
 $$\delta(Z)\ge 0.$$
The same argument shows that $$\delta(Z/A)\ge 0.$$ \qed
\df Define  for $p>0$
  $$\Z[\frac{1}{p}] =\{ \frac{m}{p^n}: m\in \Z, n\in N\}.$$
Let $A\in \SEO.$
$A$ is said to be {\bf with standard kernel} if
$$\Ker=\ker_{|A}=\{ a\in A:\ex(a)=1\}
=\left\{ \begin{array}{ll}   \pi\cdot \Z & \mbox{ if }  p=0\\
  \pi\cdot  \Z[\frac{1}{p}] & \mbox{ if }  p>0\end{array}\right.$$
for some transcendental $\pi\in D.$ \\
$A$ is said to be {\bf with compact kernel} if $\ker_{|A}$ is an algebraically
compact additive group.\\
Denote
$$\hat \Z_{(p)}=\prod_{l\neq p,\ l\mbox{ prime number}}Z_l,$$
the direct products of
additive groups of $l$-adic integers. These groups are algebraically compact
(see [F], Ch. 7).\\
$A$ is said to be {\bf with a full kernel} if for
$\ker =\ker_{|A}=\{ a\in A:\ex(a)=1\}$
the group $A/\ker$ is isomorphic to the multiplicative subgroup
of an algebraically closed field containing  all the torsion points.\\ \\
\bp \lb{kernel} For any characteristic $p$

(i) there is an $A\in \SED_p$ with standard full kernel;

(ii) there is an $A\in \SED_p$ with compact full kernel. \ep
\pf
(i) Take an additive subgroup $A=\pi\cdot\Q\subs D$ for $\Q$
the prime subfield of $D$ and $\pi$ a transcendental element in $D.$
Define  $H=A/\Ker$ for $\Ker$ the standard kernel with generator $\pi$
defined above. Then $H$ considered as multiplicative group is characterized
by the property that it is a torsion group such that any equations of the
form $x^n=h$ has for any $h$
exactly $n$ solutions when $n$ is prime to $p$ and only one when $n=p^k.$
In other words $H$ is isomorphic to  the torsion subgroup of an
algebraically closed field $R$  of characteristic $p.$  Define $\ex$ as the
canonical homomorphism $A\to R^\times$ corresponding to this isomorphism.
Since $\pi$ is transcendental, $\delta(X)=0$ for any finite $X\subs A.$\\
(ii) In a large field $D$ of characteristic $0$ consider an additive
divisible subgroup $A$ of power continuum and such
 that the generators
of $A$ as a vector space over $\Q$ are algebraically independent.
Such a group
is isomorphic to to the divisible hull of $\hat \Z(p)$ and any algebraic
dependence between elements of $A$ is a linear dependence over $\Q.$
 Fix $K\subs A$ the corresponding
subgroup isomorphic to $\hat \Z(p).$
Define $H=A/K.$  Then $H$  is a torsion group such that any equations of the
form $x^n=h$ has for any $h$
exactly $n$ solutions when $n$ is prime to $p$ and only one when $n=p^k,$
i.e. $H$ is isomorphic to  the torsion subgroup of an
algebraically closed field $R$  of characteristic $p.$  Define $\ex$ as the
canonical homomorphism $A\to R^\times$ corresponding to this isomorphism.
Since the generators of $A$ are algebraically independent, $\trd(X)=\mdim(X)$
for any finite $X\subs A.$ Thus
$\delta(X)=0$ for any finite $X\subs A.$
\qed

\bl \lb{contin}
Suppose $A \in \SED$ and $A$ is with full kernel. Then there is $D\in \EZ$
and an embedding
of $A$ into $D$ such that $A\le D$ and $\ker_{|D}=\ker_{|A}.$
\el
\pf  Choose a field $D$ of characteristic zero  and an algebraically closed
field $R$ such
that $A\subs D,$ \ $A/E\subs R^{\times}$ and  $\ldim(D/A)=\trd(R/\ex A)\ge \aleph_0.$
We want to define $\ex:D\to R^{\times}$ extending $\ex_A$ so that
 $\D\in \EZ.$

Denote $A_0=A,$ \ $\ex_0=\ex_A$ and $H_0=\ex_0(A_0).$

Proceed by induction defining $A_{\alpha},$   $H_{\alpha}$
and an endomorphism   $$\ex_{\alpha}: A_{\alpha}\to H_{\alpha}$$ by choosing

on the even steps: the first element  $a\in D\smin  A_{\alpha}$ and define
$\ex_{\alpha+1}(a)$ to be any element in $R^{\times}\smin acl( H_{\alpha}).$ Put
 $A_{\alpha+1}= A_{\alpha}+\Q\cdot a$ and extend  $\ex_{\alpha+1}$ to
 $A_{\alpha+1}$ as a group homomorphism. Put
$H_{\alpha+1}= \ex_{\alpha+1}(A_{\alpha+1}).$

on the odd steps: the first element  $h\in R^{\times}\smin  H_{\alpha}$ and define
$a$ to be any element in $D\smin \acl(A_{\alpha})$
and $\ex_{\alpha+1}(a)=h.$
 Put
 $A_{\alpha+1}= A_{\alpha}+\Q\cdot a$ and extend  $\ex_{\alpha+1}$ to
 $A_{\alpha+1}$ as a group homomorphism. Put
$H_{\alpha+1}= \ex _{\alpha+1}(A_{\alpha+1}).$

On both even and odd steps it follows from Lemma~\ref{ext} that
 $A_{\alpha+1}\in \SED$ and $A_{\alpha}\le A_{\alpha+1}.$

Also, 
$$\ker_{|A_{\alpha+1}}=\ker_{|A_{\alpha}}$$ since if
$\ex(qa+a')=1$ for some rational $q=\frac{m}{n}$ and $a'\in A_{\alpha}$
then $h^m=g^n$ for $h=\ex(a),$ $g=\ex(a')\in H_{\alpha}.$
 Since $A_{\alpha}$
is with full kernel, $ H_{\alpha}$ contains a root of degree $m$ of $g^n,$and $h\notin H_{\alpha}$ only $q=0$ is possible.
\qed \nt Let $K$ be a full kernel for some $\D$ in  $\EZ.$ Denote $\EZ(K)$ the class of structures with kernel $K.$

\ssn{Exponentially-algebraically closed structures of $\EZ(K)$}
In this section we consider  the class of structures
with a given full kernel $K.$ We extend the language $\LL_p$ by
naming all elements of $K,$ thus when we say that an algebraic variety
$V$ is defined over some $C\subs D,$ we mean the parameters in the definition
of $V$ are from the field generated by $C\cup K.$  We start by giving some basic definitions and
notations.\\ \\
\df For $C\subs D$ and algebraic varieties $V\subs D^{n},$
$W\subs R^{n}$ we say the {\bf pair $(V,W)$ is definable
over $C$ } meaning $V$ is definable over $C$ and $W$ over $\ex(span_Q(C)).$
If the varieties are irreducible over the corresponding sets, then the pair
is said to be irreducible over $C.$\\ \\
$V$ is said to be {\bf free of additive dependencies over
$C$} if  no $\bar a\in V$ generic over $C$ satisfies
$m_1\cdot a_1+\dots+m_n\cdot a_n=c$ for some $c\in \spn(C\cup K)$ and non-zero
tuple of  integers $m_1,\dots,m_n.$ $V$ is said to be {\bf aboslutely
free} of additive dependencies over $C$ if it is free over $\acl(C\cup K).$ 
$W$ is said to be
{\bf free of multiplicative dependencies over $C$} if
no $\bar b\in W$ generic over $\ex(C)$ satisfies
$ a_1^{m_1}\cdot\dots\cdot a_n^{m_n}=r$ for some $r\in \ex(\spn(C)).$
$W$ is said to be {\bf absolutely free} of multiplicative dependencies
if it is free over $\acl(\ex(C)).$\\
A pair $(V,W)$ is said to be {\bf a free  (absolutely free) pair} if
both $V$ and $W$ are free (absolutely free) of additive and
multiplicative 
dependencies correspondingly.  \\ \\
Let $W\subs R^n$ be an algebraic variety defined and irreducible over
some $\ex(C)$ for some $C=span_Q(C\cup K)\subs D.$
With any such $W$  we associate a sequence $\{ W^{\frac{1}{l}} :l\in \N\}$ of
algebraic varieties which are definable and irreducible over $\ex(C)$ and satisfy
the following:\\
$W^{1}=W$ and for any $l,m\in \N$ the mapping
$$[m] :\la y_1,\dots y_n\ra \mapsto \la y_1^m,\dots y_n^m\ra$$
maps $W^{\frac{1}{lm}}$ onto $W^{\frac{1}{l}}.$ Such a sequence is said to
be {\bf the sequence   associated with $W$} over $C.$\\
Also with any $\la y_1,\dots y_n\ra\in W$ as above we associate a sequence
$$\{ \la y_1,\dots y_n\ra^{\frac{1}{l}}: l\in \N\}$$
such that for any $l,m\in \N$ the mapping
$$[m]:\la y_1,\dots y_n\ra \mapsto \la y_1^m,\dots y_n^m\ra$$
maps $\la y_1,\dots y_n\ra^{\frac{1}{lm}}$
onto $\la y_1,\dots y_n\ra^{\frac{1}{l}}.$
Such a sequence is said to be
{\bf associated} with $\la y_1,\dots y_n\ra.$\\ \\
Let $(V,W),$ $(V',W')$ be  pairs over $C,$ $(V,W)$ irreducible,
$V'\subsetneqq V,$ $W'\subsetneqq W,$
$\{ W^{\frac{1}{l}}:l\in \N\}$  a sequence associated with $W.$ Then the triple
$$\tau=(V\smin V',W\smin W',\{ W^{\frac{1}{l}}:l\in \N\})$$ is said to be {\bf
an almost finite [quantifier-free]  type over} $C.$
A pair of the form $(V\smin V',W\smin W')$ is said to be a {\bf finite type}
over $C.$

A tuple $\la a_1,\dots a_n\ra\in D^n$ is said to be realising  type $\tau$
if $\la a_1,\dots a_n\ra\in V\smin V',$
$\ex\la a_1,\dots a_n\ra\in W\smin W'$ and
$\ex\frac{1}{l}\la a_1,\dots a_n\ra\in W^{\frac{1}{l}}$ for all
$l\in \N.$  \\ \\
A pair $(V,W)$ is said to be a {\bf normal pair over $C$ } if in some extensions of
the fields  there are  $\la a_1,\dots ,a_{n}\ra\in V$ and $\la b_1,\dots ,b_{n}\ra\in W$
 such that for any
$k\le n$ independent integer vectors $m_i=\la m_{i,1},\dots
m_{i,n}\ra,$\ $i=1,\dots, k,$\ and
 $$a'_i=m_{i,1}a_1+\dots+ m_{i,n}a_n,\ \ \
b'_i=b_1^{m_{i,1}}\cdot\dots \cdot b_n^{m_{i,n}}$$
it holds
$$\ldim (\la a'_{1},\dots, a'_{k}\ra/C) +
\trd(\la b'_{1},\dots, b'_{k}\ra /\ex (C))
\ge k.$$
Equivalently, for the varieties
$$V'_{1,\dots, k}={\rm locus}_C(a'_1,\dots,a'_k)\mbox{ and }
W'_{1,\dots, k}={\rm locus}_{\ex(C)}(b'_1,\dots,b'_k)$$
it must hold
$$\dim V'_{1,\dots, k}+\dim W'_{1,\dots, k}\ge k.$$

Notice that in case $$\la a'_1,\dots, a'_k\ra=\la a_{i_1},\dots,
a_{i_k}\ra
\mbox{ and }\la b'_1,\dots, b'_k\ra=\la b_{i_1},\dots, b_{i_k}\ra$$
the varieties $ V'_{1,\dots, k}$ and $W'_{1,\dots, k}$ are just the
projections
$ V_{i_1,\dots, i_k}$ and $W_{i_1,\dots, i_k}$ of the initial
varieties on the $({i_1,\dots, i_k})$-subspaces.\\

\rem Given a non-degenerate integer matrix $M=\{ m_i:i=1,\dots,m_n\},$
we consider the mapping $$[M]:\bar a\bar b\mapsto \bar a'\bar b'$$
determined by the $n$ integer vectors of the matrix as in the
definition above. This is a linear transformation on $D^n$ and a
rational transformation on $R^n.$ We will call  the mapping a {\bf
linearly induced mapping} $[M].$\\

\df A structure $\D$ in $\EZ(K)$ is said to be {\bf $\EZ(K)$-exponentially-algebraically closed}
(e.a.c.) if for any $\D'\in \EZ(K),$ such that $\D\le \D',$  any  almost finite
quantifier-free type
over $D$ which is realised in $\D'$   has a realisation in $\D.$\\
The class of $\EZ(K)$-exponentially-algebraically closed structures is denoted
$\EC(K).$ We omit $K$ when the kernel is fixed.\\ \\
We continue by studying properties of associated sequences.
A sequence associated with  $\bar y$ is not uniquely determined;
for $\bar y^{\frac{1}{l}}$ there are  $l^n$ possible values. Evidentely,
one can get all the values  multiplying a value $\bar y^{\frac{1}{l}}$ by
all $\bar \xi=\la \xi_1,\dots,\xi_n\ra,$ where $\xi_i$'s are roots
of unity of degree $l.$ It follows
\bl \lb{n1} Let
$\bar y\in W$ and $\{ \bar y^{\frac{1}{l}}\in l\in \N\},$
$\{ W^{\frac{1}{l}} :l\in \N\}$ be sequences associated with $\bar y$ and
$W$ correspondingly. Then there is a sequence $\{ \bar \xi(l):l\in \N\}$
of roots of unity of power $l,$ such that
$$\bar \xi(l)\cdot \bar y^{\frac{1}{l}}\in W^{\frac{1}{l}}$$ for all
$l\in \N$ and $\xi(l\cdot m)^m=\xi(l)$ for all $m\in \N.$\el
\bl \lb{n2} Assume $\D\in \SED$ is with a compact kernel,
$W$ a nonempty algebraic subvariety of $R^n$   and
$\{ W^{\frac{1}{l}} :l\in \N\}$ a sequence associated with $W.$
Then there is $\bar x\in D^n$ such that
$$\ex({\frac{1}{l}}\cdot \bar x)\in W^{\frac{1}{l}}$$ for all $l\in \N.$\el
\pf Take any $\bar y\in W$ and $\bar v\in D^n$ such that $\ex(\bar v)=\bar y.$
Evidentely, $\{ \ex({\frac{1}{l}}\cdot \bar v):l\in \N\}$ is a sequence
associated with $\bar y.$
By Lemma~\ref{n1} there is a sequence $\{ \bar \xi(l):l\in \N\}$ such that
$\bar \xi(l)\cdot \ex({\frac{1}{l}}\cdot \bar v)\in W^{\frac{1}{l}}$ for all
$l\in \N.$ Let $\bar \alpha(l)\in D^n$ be such that
$\ex(\bar \alpha(l))=\bar \xi(l)$ for all $l\in \N.$
Consider a system of equations in the additive language in variables
$\bar \chi,$ $\bar z(l)$ ($l\in \N$) running in $\ker^n:$
$$\{ \bar \chi=l\cdot \bar \alpha(l)+l\cdot \bar z(l):l\in \N\}.$$
Any finite subsystem
$$\{ \bar \chi=l\cdot \bar \alpha(l)+l\cdot \bar z(l):l\le m\}$$
of the system has a solution $\bar \chi=m!\bar \alpha(m!).$ Indeed,
$$\ex({\frac{1}{l}}\cdot \bar \chi)=\bar \xi(m!)^{\frac{m!}{l}}=\bar \xi(l)=\ex(\alpha(l)),$$
hence $\bar \chi-l\cdot \bar \alpha (l)\in l\cdot\ker.$ Thus the whole system
has a solution $\bar \chi\in \ker^n.$ Put $\bar x=\bar \chi+\bar v.$\qed \\
 If the number of the varieties
 $W^{\frac{1}{l}}$ for a given $W$ grows with $l,$ then by the symmetry
coming from the action of Galois group (multiplication by roots of unity of
degree  $l$), it follows that the number of the sequences associated with $W$
 is  $2^{\aleph_0}.$ This is the case when $W$ is not free of multiplicative
 dependence. Indeed, if $y_1^{m_{i,1}}\cdot\dots\cdot y_1^{m_{i,n}}=c$
holds on $W$ then for $W^{\frac{1}{l}}$ holds
$y^{m_{i,1}}\cdot\dots\cdot y^{m_{i,n}}=c^{\frac{1}{l}}$ with exactly $l$
different choices for $c^{\frac{1}{l}}.$ Nevertheless, in case $W$ is free
of multiplicative dependencies  F.Voloch answered our question by proving

\bp \lb{Voloch}[F.Voloch]  Let $W$ be an   irredicible algebraic variety
such that $W$ is not  contained in a coset of any proper torus.   Then
there exist a number $c$ such that the number of irreducible components of
$W^{\frac{1}{l}}$ is bounded by $c$ independently on $l.$\ep
\pf First we reduce the proof to the case when $W$ is a curve by cutting $W$
by a generic linear subspace of appropriate dimension.
Denote $K=R(W)$ the field of rational functions on the curve $W$
defined over  $R.$ Let $x_1,\dots, x_n\in K^{\times}$ be the coordinate
functions, which are  multiplicatively independent over
$R^{\times}$ by assumptions. \\
Claim. Let $H$ be the group generated by
$x_1,\dots, x_n$ in $K^*.$ Then there is $d>0,$ such that
 $\# H/H\cap(K^*)^l\ge dl^n$ for any $l\ge 1.$\\
Indeed, by hypothesis $x_1,\dots, x_n$ generate a free abelian group of
rank  $n$  in $K^*/R^*$ isomorphic to $H,$ and so $\# H/H^l=l^n.$
So we need to show that $\# (H\cap K^l/H^l)$ is bounded independent of $l.$
Let ${\cal D}$ be the group of divisors of $K/R$ and ${\cal D}_S$ the
subgroup generated by $S=\{ v| \exists i \ \ v(x_i)\neq 0\}.$ So ${\cal D}_S$
 is a f.g. abelian group and $H$ embeds in ${\cal D}_S.$ Let
$$G={\cal D}_S\cap (H\otimes \Q)\subs {\cal D}_S\otimes \Q$$ which is
also a lattice and $H\subs G$ with $G/H$ finite, $\# G/H=c.$ Now, if
$x\in H\cap (K^*)^l$ then $x=y^l,$ $y\in K^*.$  But $(y)\in {\cal D}_S$
and also $H\otimes \Q,$ so $(y)\in G,$ so $y^c\in H,$ so $x^c=y^{cl}\in H^l,$
hence $H\cap K^l/H^l$ is an abelian group of exponent $c$ generated by
$\le n$ elements, so $\# H\cap K^l/H^l\le c^n$ and $d= 1/c^n.$ Claim proved.

Let $W^{\frac{1}{l}}$ be an irreducible curve such that
$$[l](W^{\frac{1}{l}})=W.$$  Then for the field
$R(W^{\frac{1}{l}})$ of rational functions of
the curve by Galois theory  $$|R(W^{\frac{1}{l}}):K|= \# H/H\cap(K^*)^l.$$
The latter is bounded from below by the estimate $dl^n,$ by the claim.
On the other hand the mapping $[l]$ on the full inverse $[l]^{-1}(W)$
by Bezout theorem is  of degree $l^n.$
It follows that the number of irreducible components in the inverse
image is at most $1/d=c.$   \qed

\bc \lb{c2} If $W$ is an algebraic variety absolutely free of multiplicative dependencies
then any irreducible component of $W$ satisfies assumptions of
Proposition~\ref{Voloch} and thus there is an $l_0\in \N$ such that sequences
$\{ W^{\frac{1}{l}}:l\in \N\}$ and
$\{ \bar y^{\frac{1}{l}}:l\in \N\}$ satisfy
$$\bar y^{\frac{1}{l}}\in W^{\frac{1}{l}}\mbox{ for all } \ l\in \N\mbox{ iff }
\bar y^{\frac{1}{l}}\in W^{\frac{1}{l}}\mbox{ for all }l\le l_0.$$
In other words, the sequence $\{ W^{\frac{1}{l}}:i\in \N\}$ is
determined by its cut of length $l_0.$ \ec

\bl \lb{1/l} Let $W\subs R^{n+m}$ be an  algebraic variety definable and
irreducible over some subfield $F$ and absolutely free of multiplicative dependencies. Let
$A$ be a torsion-free divisible subgroup of the additive group of $D$
of rank $m,$  $\bar a\in A^m,$ $W_a=W(\ex(\bar a)).$ Then over
$\ex(A)\cdot F$ any associated sequence
$\{ W_a^{\frac{1}{l}}:l\in \N\}$ is determined by its finite cut.\el
\pf    We may assume $\bar a$ generates $A$ over $\Q.$ 

Let  $\{ \bar b^{\frac{1}{l}}:l\in \N\}$
be a realisation of an associated sequence $\{ W_a^{\frac{1}{l}}:l\in \N\}.$
Then the algebraic locus of $\bar b^{\frac{1}{l}}$ for an $l\in \N$ over
$\ex(A)$ is of the form $W'(\ex(\bar a'))$ for some $\bar a'\in A^m,$
 $W'\subs R^{n+m}$ over $F.$
Since $\bar a$ generates $A,$ for some $k\in \N$ \ $\ex(\bar a')$ is a rational
combination of $\ex({\frac{1}{lk}}\bar a).$ Hence we may assume
$\bar a'={\frac{1}{lk}}\bar a$ in the expression above. On the other hand
the locus of $(n+m)$-tuple $\bar b^{\frac{1}{lk}}\ex({\frac{1}{lk}}\bar a)$
over $F$ is given by a member $W^{\frac{1}{lk}}$ of a sequence associated with
$W.$ In other words $W'(\ex(\bar a'))$ can be represented as
$(W^{\frac{1}{lk}}(\ex({\frac{1}{lk}}\bar a))^k$ for some $k\in N.$
Suppose $l\ge l_0.$ Then $W^{\frac{1}{lk}}$ is determined uniquely once
$W^{\frac{1}{l}}$ and $k$ is given. It follows
$$(W^{\frac{1}{lk}}(\ex({\frac{1}{lk}}\bar a))^k=
W^{\frac{1}{l}}(\ex({\frac{1}{l}}\bar a).$$
Finally notice that the sequence
$\{ W^{\frac{1}{l}}(\ex({\frac{1}{l}}\bar a):l\in \N\}$ is determined by its
$l_0$-cut. \qed

\nt For subsets $A\subs B\subs D$ with $\D\in \EZ(K)$ the notion
$A\le B$ is applied in the same sense as for substructures.\\ \\
Now we come to study the rank notion $\dd_D$ for $\D$ e.a.c..

\bl \lb{part} For $\D\in \EC$, given finite $A\subs D$
and $\D'\in \EZ(K)$ such that $D\le D',$
$$\partial_D(A)=\partial_{D'}(A).$$\el
\pf Immediate from the definitions.\qed\\
By Lemma above we may omit $D$ when writing $\partial_D.$
\bl \lb{star}
For any  finite $A\subs D$ and any $a,b\in D$

(i) $\partial(A)\le \partial(aA)\le \partial(A)+1;$

(ii) $\partial(abA)=\partial(aA)=\partial(A) \mbox{ implies } \partial(bA)=\partial(A);$

(iii) $\partial(abA)=\partial(aA)\& \partial(A)< \partial(bA)\mbox{ implies } \partial(abA)=\partial(bA);$

(iv) $\partial(aA)=\partial(A)=\partial(bA) \mbox{ implies } \partial(abA)=\partial(A);$

(v) $\partial(aA)=\partial(A) \mbox{ implies } \partial(bA)=\partial(abA).$\el
\pf (i) follows immediately from the definitions of $\delta$ and $d.$
(ii) and (iii) are immediate from (i).

(iv) Let $B'\sups aA,$  $B''\sups bA$ be such that $\delta (B')=\partial(aA)$
and $\delta (B'')=\partial(bA).$ Denote $B=B'\cap B''.$

Notice that $\delta (B'\cup B'')\le \delta (B'').$ Indeed by
Lemma~\ref{delta}(iii)\\
$\ \ \ \ \delta (B'\cup B'')=\delta (B'/B'')+\delta (B'')=$\\
$=[ \ldim(B'/B'')+\trd(\ex(B')/\ex(B''))-\mdim(B'/B'')]+\delta (B'').$\\
By the modularity of linear dimension $\mdim(B'/B'')=\mdim(B'/B).$
Also, by properties of algebraic dependence \\
$\trd(B'/B'')\le \trd(B'/B)$ and $\trd(\ex(B')/\ex(B''))\le \trd(\ex(B')/\ex(B))$ Hence $\delta (B'/B'') \le \delta (B'/B).$
The latter is less or equal to zero by the choice of $A,$ $B'$ and $B.$

Now, since $abA\subs B'\cup B''$ and
$\delta (B'\cup B'')\le \delta (B'')=\partial(A),$ we have $\partial(abA)=\partial(A).$

(v) is immediate from (iv).\qed
\bc \lb{cl} The operator $$A\mapsto \cl(A)=\{ b: \partial(Ab)=\partial(A)\}$$
in $\D$ is a closure operator.
$\cl(A)$ is said to be the  $\partial$-closure of $A.$\ec
\df Let $C\subs D$ and $\bar a\in D^n$ be given. The
{\bf \ex-locus of $\bar a$
over $C$} is given as $$(V,W, \{ W^{\frac{1}{l}}:l\in \N\}),$$ where
 $V\subs D^n$ is the minimal algebraic variety over $C\cup K$ containing $\bar a$
and $W^{\frac{1}{l}}\subs R^n$ the minimal algebraic variety over
$\ex(\spn (C))$
containing $\ex (\frac{1}{l}\bar a).$ We will  call $(V,W)$ the [incomplete]
 \ex-locus.
\bl \lb{normal} Let $C, A\in \SED,$  $C\le A,$
 and let $\bar a$ be a linearly independent over $C$
string of elements of $A.$ 
Then the incomplete \ex-locus  $(V,W)$ of $\bar a$ over $C$ is
normal, $V$ is  free of additive dependencies, and absolutely free if
$C$ is an algebraically closed subfield of $D.$ 
If  $C$ is with full kernel and $\ex(C)$ is a subfield of $R,$ then $W$ is free of
multiplicative dependencies. If $\ex(C)$ is algebraically closed, then
$W$ is absolutely free.\el
\pf Let $\la b_1,\dots,b_{n}\ra$ in the definition of normality be
$\ex(\bar a).$ Then the inequalities in the definition of normality follow from
$C\le A.$

An additive dependence for $V$ would mean by definition of $V$ a
linear dependence of $\bar a$ over $C,$ which does not hold.
 
A multiplicative dependence for $W$ is equivalent to $\ex(\bar a)$
being multiplicatively dependent over the subfield generated  by
$\ex(C),$ which is equivalent under the assumptions to $\bar a$ being
linearly dependent over $C.$  
\qed

\bl \lb{almost} Suppose $\D\in \EC,$  $\D\le \D'.$ Let $\bar a$ be a finite
string in $D'$ $\tau=(V\smin V',W\smin W',\{ W^{\frac{1}{l}}:l\in
\N\})$ and
$(V,W,\{ W^{\frac{1}{l}}:l\in \N\})$  the ex-locus of $\bar a$ over $D.$ Then
there is a realisation of $\tau$ in $\D.$\el
\pf We may assume $\bar a$ is linearly independent over $D.$ Then by
Lemma~\ref{normal} $(V,W)$ is  normal and absolutely free pair,
since $D$ and $R$ are algebraically closed fields. By
Corollary~\ref{c2} $\tau$ is equivalent to a finite type. By the
assumption for $\D,$  $\tau$ is realised in $\D.$\qed 

 \bl \lb{ec1}  Let $K\subs C\subs D,$ 
$\tau=(V\smin V',W\smin W',\{ W^{\frac{1}{l}}:l\in \N\})$ 
be an almost finite type over $C$ and assume the
pair $(V,W)$ is normal over $C,$ $V$ is free of additive dependencies
and $W$ is absolutely free of multiplicative dependencies over $C.$
Then there is $\bar a$ in $\D$ realising
$\tau.$
Moreover in some extension $\D'\ge \D$ \
$\bar a$ can be chosen generic in $V$
(of maximal  tr.d.) over $C$ and $\ex(\bar a)$
generic in $W$ over $\ex(C).$\el
\pf Let $\D'$ be a normal extension of $\D$ with $D'$ and $R'$ of
infinite transcendence degree.
Applying a linearly induced transformation we may assume that for any
$\bar a= \la a_1,\dots,a_n\ra\in V$ in $D',$
generic   over $C,$   $ a_1,\dots,a_k$ are linearly independent over
$\acl(C)$ and  $a_{k+1},\dots,a_n\in \acl(C).$ And since $V$ is free
over $C$   $a_1,\dots,a_n$ are linearly independent over $C.$

Denote $W_{k+1,\dots,n}$ the variety induced by $W$ on 
$\{ k+1,\dots,n\}$-coordinates. It follows from normality that 
$\dim W_{k+1,\dots,n}=n-k$ and thus by letting
$b_i^{\frac{1}{l}}=\ex(\frac{1}{l}a_i)$ for each $k<i\le n$ and each
$l\in \N$ we get      
$\la b_{k+1},\dots,b_n\ra^{\frac{1}{l}}\in
W_{{k+1},\dots,n}^{\frac{1}{l}}$ generically. 
Now extend $\la b_{k+1},\dots,b_n\ra^{\frac{1}{l}}$ to 
$\{ \la b_{1},\dots,b_n\ra^{\frac{1}{l}}\in  W^{\frac{1}{l}}$ for each
$l\in \N.$

By absolute freeness,  
$b_1,\dots,b_n$ are multiplicatively independent over $\acl(\ex(C)).$\\

Extend $\ex$ to $A=D+\spn(a_1,\dots a_n)=D+\spn( a_1,\dots,a_k)$ as:
$$\ex(\sum_{i\le k}\frac{m_i}{l}a_i+d)=\prod_{i\le k}({b_i}^{\frac{1}{l}})^{m_i}\cdot \ex(d)$$
for any integers $m_i,$ $l\neq 0$ and element $d\in D.$
The definition is consistent since  $a_1,\dots,a_k$ are linearly
independent over $\acl(C)$ and hence over $D.$
So the formula
defines a homomorphism. The kernel of the homomorphism coincides with that of
$\ex$ on $D,$ since $\prod_{i\le k}({b_i}^{\frac{1}{l}})^{m_i}\cdot
r=1$ for some $m_i$ and $r\in R$ implies by genericity $r\in
\acl(\ex(C))$ and hence by the assumptions $m_1=\dots=m_k=0$ and $r=1.$
 Thus $A\in \SEO$ and is with full kernel.
Notice that by
the normality 
 for any
independent integer vectors $m_i=\la m_{i,1},\dots
m_{i,n}\ra,$\ $i=1,\dots, k,$
  \ \ it holds
$\delta(m_1\bar a,\dots, m_k\bar a/D)\ge 0.$

Thus $D\subs A$ satisfy the assumptions of Lemma~\ref{ext} and hence
$A\in \SED,$ $D\le A.$
By the choice $\bar a$ realises $\tau.$
Since $\D\in \EC$ by Lemma~\ref{almost} there is a  realisation of the type  in $\D.$
\qed

\bl \lb{3.7+}  Let $K\subs C\subs D,$ 
 and assume the
pair $(V,W)$ is free and normal over $C.$ Then for any $V'\subsetneqq
 V$ and $W'\subsetneqq W$ over $C$
there is $\bar a\in V\smin V'$ living in $\D$ such that 
$\ex (\bar a)\in W\smin W'.$ 
Moreover in some extension $\D'\ge \D$ \
$\bar a$ can be chosen generic in $V$
 over $C$ and $\ex(\bar a)$
generic in $W$ over $\ex(C).$

If the kernel is compact then for any sequence $\{ W^{\frac{1}{l}}:
l\in \N\}$ we can find a realisation $\bar a$ of type $(V\smin V',W\smin W',\{ W^{\frac{1}{l}}:l\in \N\}).$ 
\el
\pf We may again assume that for any generic  $\bar a\in V$ in some
 extension \ 
$a_1,\dots,a_k$ are linearly independent over $\acl(C)$ and
 $a_{k+1},\dots,a_n\in \acl(C)\subs D.$
Also by transformations we may assume that for any generic $\bar b$ \ $b_1,\dots,b_s$ are
multiplicatively independent over $R$ and $ b_{s+1},\dots,b_k\in
 \acl(\ex(C))\subs R$ for some $s\le k.$ 
It follows now from normality
 that the algebraic type of $ a_{s+1},\dots,a_k$ is characterised by
 the property that the elements are algebraically independent over 
$C.$

Denote $W_{s+1,\dots,k}$ the variety induced by $W$ on 
$\{ s+1,\dots,k\}$-coordinates.  
Now, choose first a generic $\la b_{s+1},\dots,b_k\ra\in
W_{{s+1},\dots,_k}$ over $\ex(C).$ By the above noted
 $ b_{s+1},\dots,b_k\in R.$ In case the kernel is compact consider
also an associated sequence 
$$\{ \la b_{s+1},\dots,b_k\ra^{\frac{1}{l}}\in  W^{\frac{1}{l}}_{s+1,\dots,k}:l\in \N\}.$$
Now
choose for
$s<i\le k$ \ $a_i\in D$
so that $\ex(a_i)=b_i$ and, in the case of compact kernel, by
Lemma~\ref{n2}, $\ex({\frac{1}{l}}a_i)=b_i^{\frac{1}{l}}.$
In the other case just denote   
$\ex(\frac{1}{l}a_i)=b_i^{\frac{1}{l}}.$ 

Choose now $\la a_{k+1},\dots,a_n\ra \in V_{k+1,\dots,n},$ the variety
induced by $V$ on the last $n-k$ coordinates. By the above noted these
are in $R$ so we can define for $k<i\le n$ \ \
$b_i^{\frac{1}{l}}=\ex({\frac{1}{l}}a_i).$ Then $\la
b_{s+1},\dots,b_n\ra$ will satisfy $W_{s+1,\dots,n}$ since the last
$n-k$ coordinates are algebraically independent over the preceding
ones. In case the kernel is compact, by the choice of $ \la
b_{s+1},\dots,b_k\ra^{\frac{1}{l}}$ this is also true for
$W_{{s+1},\dots,n}^{\frac{1}{l}}.$

Thus $\la a_{s+1},\dots,a_n\ra$ realises the type
$\tau_{{s+1},\dots,n}$ induced by $\tau$ on the last $n-s$
coordinates.

Finally notice that the type $\tau(a_{s+1},\dots,a_n)$ over $C\cup
\{ a_{s+1},\dots,a_n \},$ corresponding to the first $s$ coordinates
in $\tau,$ satisfies the assumptions of Lemma~\ref{ec1}. Thus it has a
realisation $\la a_1,\dots, a_s\ra$ in $\D.$ Which completes the
construction of $\la a_1,\dots,a_n\ra.$ \qed 

Combined with Lemma~\ref{normal} we thus get
\bc \lb{c3.7+} 
(i) The set of incomplete ex-loci of $\bar a$ over $C\le
\D\in \EC$ for $\bar a$ living in $\D$ does not depend on $\D.$

(ii) If the kernel $K$ is compact then the set of (complete) ex-loci of $\bar a$ over $C\le
\D\in \EC$ for $\bar a$ living in $\D$ does not depend on $\D.$
\ec

\bt \lb{t1} A structure $\D\in \EZ(K)$ is in $\EC(K)$ iff for any irreducible
normal  free pair
$(V,W)$ over
$\D$  there is a
realisation of the type $(V,W)$ in $\D.$\et
First we prove
\bl \lb{V'}
For any irreducible free normal pair $(V,W)$ ($V\subs D^n, W\subs R^n$)
and $V'\subs V,$ $W'\subs W$ there is a free normal pair
$(V^*,W^*)$ ($V^*\subs D^{n+m}, W\subs R^{n+m}$)
such that $\la a_1,\dots,a_n\ra\in D^n$ realises $(V\smin V',W\smin W')$ iff
there is $\la a_{n+1},\dots, a_{n+m}\ra \in D^m$ such that
$\la a_1,\dots,a_n,a_{n+1},\dots,a_{n+m}\ra$
realises $(V^*,W^*).$ \el
Let $f(x_1,\dots,x_n)$
be a polynomial in ${\rm Ann} V'\smin {\rm Ann} V.$
Choose a natural number $k$
such that there is no rational function $r$ such that $r^k\equiv f$ on  $V.$
Add new variable $x_{n+1}$ together with new identity
$f(x_1,\dots,x_n)\cdot x_{n+1}^k=1,$ which
defines together with   ${\rm Ann} V$ new algebraic variety $V^f.$
Notice that the projection of $V^f$ onto the first $n$ coordinates is
$V\smin V_f$ where $V_f$ is the subvariety of $V$
defined by $f=0.$ Also, $V^*$ is free of additive dependencies since
otherwise  $\sqrt[k]{f}=x_{n+1}$ is a rational function on $V.$
Using this method finitely many times we come to $V_1^*\subs D^{n+l}$ with the projection
$V\smin V'.$ 

Let $W_1$ be just $W\times D^{l}.$ By construction
$(V_1^*,W_1^*)$ satisfies the statement of the lemma for $(V^*,W^*)$
in case $W'=\emptyset.$\\

Assume now $W'\neq \emptyset.$  Let $g(y_1,\dots,y_n)$  be a polynomial in ${\rm Ann} W'\smin {\rm Ann} W.$
We may assume that there is no integer $k$
and $m_1,\dots, m_n$ not all zero such that
$$g(y_1,\dots,y_n)^k\cdot y_1^{m_1}\cdot\dots\cdot y_n^{m_n}$$  is
constant on $W$ since otherwise $g(y_1,\dots,y_n)$ is nonvanishing on $W.$ Add new variable $y_{n+1}$ together with new identity
$g(y_1,\dots,y_n)\cdot y_{n+1}=1,$ and the resulting variety $W_1^g$
projects onto 
the first $n$ coordinates as $W\smin W_g$ where $W_g$  is the subvariety of
$W$ defined by $g=0.$ By our assumption $W_1^g$ is free of multiplicative
dependencies and we can come to $W^*$ as above.

The normality of the pair is evident from the construction.\qed 

\pf of Theorem. The
left-to-right implication follows
 from Lemma~\ref{ec1}.
To get the inverse assume $\bar a$ is in some $\D'\ge \D$ and we need to
realise an almost finite type
$(V\smin V',W\smin W',\{ W^{\frac{1}{l}}:l\in \N\})$
where $\tau=(V,W,\{ W^{\frac{1}{l}}:l\in \N\})$ is the ex-locus of
$\bar a$ over 
$D.$ It is enough to solve the problem for a $\Q$-linear basis $\bar a_0$ of
$\bar a$ over $D,$ so we may assume $\bar a$ is $\Q$-linearly independent
over $D.$ Thus $(V,W)$
is a free normal pair, and so we may assume $\tau$ is a finite type.  
By Lemma~\ref{V'} we reduce the type
$(V\smin V',W\smin W')$ to a type of the form
$(V,W)$ and $(V,W)$ normal and free.
By the assumptions of Theorem the type is realised in $\D.$ \qed
\bc The structure $(\C,\exp,\C)$ on complex numbers is in $\EC$ iff it
satisfies the weak Schanuel conjecture and  for any normal free pair
$(V,W)$ over $\C$ there is a realisation of the type $(V,W)$ in $\C.$\ec

\ssn{Atomic and prime structures}
In this section $\D$ is always exponentially-algebraically closed with
a fixed kernel, the
sort $R$ is of characteristic zero. We consider two cases: the
standard kernel $\pi \Z$ and the canonical compact kernel $\pi \hat
\Z.$  
\df Let $C_0\subs C_1\subs C\le D$ and $V, V'\subs D^n$   algebraic varieties
over $C_0,$ \   $W, W'\subs R^n$ varieties  over $\ex(C_0),$ \ \
$V$ and $W$  irreducible over corresponding sets and an associated
sequence $\{ W^{\frac{1}{l}}:l\in \N\}$ defined over $C_1.$
The
type $\tau=(V\smin V',W\smin W',\{ W^{\frac{1}{l}}:l\in \N\})$ is said to be a
{\bf q-atom over $C$ weakly defined over $C_0,$ defined over $C_1,$}
 if there is a realisation  $\bar a$ of $\tau$ in some
$\D'\ge \D$
and for every  $\bar a$  realising $\tau$ in a
$\D'\ge \D,$
$\bar a$ {\bf is generic over $C$ realisation of }$\tau,$
i.e.  $\bar a$ generic in $V$ over $C$ and $\ex(\bar a)$ generic in $W$
over $\ex(C).$\\

$C_0$  is said to be weakly isolating $\bar a$ over $C$ and $C_1$
isolating $\bar a$ over $C.$

 Any realisation $\bar a$ in $\D$ of a type
as above is said to be {\bf q-atomic over }$C.$\\ \\
\rem It follows immediately from definitions that for $\bar a$ and
$\tau$  as above q-atomic over $C,$ \ $\ex(\frac{1}{l}\bar a)$ is also
generic in  $W^{\frac{1}{l}}$ for any $l\in \N.$\\

{\bf Remark} It follows easily from definitions that 
if $\bar a$ and $\bar a'$ satisfy same q-atomic type over
$C,$ then the substructures $\spn(C\cup |\bar a|)$ and $\spn(C\cup
|\bar a'|)$ are isomorphic.\\ \\

\nt For $C\subs D, $ $\D\in \EC,$ denote $\ecl_D (C)$ the field obtained by the following
process in $\omega$ steps carried on inside $\D:$

$C^0=\acl(C),$ $C^{(2n+1)}=\ln(\acl(\ex(C^{(2n)}))),$ $C^{(2n+2)}=\acl(C^{(2n+1)}),$
$$\ecl_D(C)=\bigcup_{n<\omega} C^{(n)}.$$
We omit the subscript when $\D$ is fixed.\\

{\bf Remark}. If $X\subs C^{(m)},$ then by the definition
$$\delta(X/ C^{(m-1)})\le 0$$
(here $C^{(-1)}=C$). If $C^{(m-1)}\le D,$  then the equality holds.
\bl \lb{ecl}

(i) For any $n\in \N,$  $C\le D,$ $X\subs C^{(n)}$ there is
$X\subs X'\subs C^{(n)},$ such that $\delta(X'/C)=0;$

(ii) If $C\le D,$ then $C^{(n)}\le D$ for all $n\in \N$ and $\ecl(C)\le D;$

(iii) If $X\subs \ecl(C)$ and $C\le CX\le D,$ then $\delta(X/C)=0.$
\el
\pf (i) The statement is obvious for $n=0.$ So, we assume $n>0.$
Let $X\subs C^{(n)}$ and $Y\subs C^{(n-1)}$ be such that
$$\delta(X/ C^{(n-1)})=\delta(X/CY)=0.$$
By induction hypothesis we may assume $\delta(Y/C)=0.$
Since
$$\delta(XY/C)=\delta(X/CY)+\delta(Y/C)$$
we get the required by putting $X'=XY.$\\
(ii) Let $X\subs D$ be finite. We want to show that
$\delta(X/ C^{(n)})\ge 0.$ By definition $\delta(X/ C^{(n-1)})=\delta(X/CY)$
for some finite $Y\subs  C^{(n)}.$ By (i) we may choose $Y$ so that
$\delta(Y/ C)=0.$ Then
$\delta(X/CY)=\delta(XY/C)\ge 0,$ since $C\le D.$\\
(iii) follows from (i). \qed
\bl \lb{401}
Let $C\le D,$ $|\bar a|\subs \ecl(C)$
 and $C\bar a\le D.$ Then
$\delta(\bar a/C)=0$ and  $\bar a$ is
q-atomic over $C.$\el
\pf $\delta(\bar a/C)=0$ by Lemma~\ref{ecl}(iii).
To prove the rest of Lemma we first prove\\
{\bf Claim.}

(i) If $|\bar a|\subs C^{(2n+2)},$ then there is
$\bar b$ in $C^{(2n+1)}\cap \spn_Q(C\bar a),$ such that $\bar a\in \acl(C\bar b)$
and $\delta(\bar b/C)=0.$

(ii) If $|\bar a| \subs C^{(2n+1)},$ then there is
$\bar b$ in $C^{(2n)}\cap \spn_Q(C\bar a),$ such that
$\ex(\bar a)\in \acl(\ex(C\bar b))$
and $\delta(\bar b/C)=0.$\\
\pf of (i). Choose  finite
$\bar b$  to be a $\Q$-basis of $C^{(2n+1)}\cap \spn_Q(C\bar a)$
Then $\mdim(\bar a/C\bar b)=\mdim(\bar a/C^{(2n+1)}),$ since if a $\Q$-linear
combination $u$ of $\bar a$ belongs to $C^{(2n+1)},$ then
$u\in C^{(2n+1)}\cap \spn_Q(C\bar a).$
It follows that $\delta(\bar a/C\bar b)=\trd(\bar a/C\bar b)\ge 0,$ since
$\trd(\ex(\bar a)/\ex(C\bar b))\le \mdim (\bar a/C\bar b)$ and\\
$\trd(\ex(\bar a)/\ex(C\bar b))\ge \trd(\ex(\bar a)/\ex(C^{(2n+1)}))=
\mdim(\bar a/C^{(2n+1)}) =\mdim(\bar a/C\bar b).$ \\
On the other hand, $\delta(\bar a/C\bar b)=\delta(\bar a\bar b/C)-
\delta(\bar b/C)=0=\delta(\bar b/C), $ since
$\delta(\bar a\bar b/C)=\delta(\bar a/C)=0$ and
$\delta(\bar b/C)\ge 0$ by $C\le D.$ Hence $\trd(\bar a/C\bar b)=0,$ which
means $|\bar a|\subs \acl(C\bar b).$\\
\pf of (ii) is symmetric.    \\
Now we continue the proof of Lemma. By induction on $m$ we prove that
if $|\bar a|\subs C^{(m)},$ then $\bar a$ is q-atomic over $C.$
For $m=0$ \ $\bar a\in \acl(C)$ and the statement is evident.
If $|\bar a|\subs C^{(m)}$ for $m>1,$ then by the Claim there is
$\bar b$ in $C^{(m-1)}\cap \spn_Q(C\bar a)$ such that $\bar a\in \acl(C\bar b)$
or $\ex(\bar a)\in \acl(\ex(C\bar b)).$  Let
$\bar a\in \acl(C\bar b)\subs (C\bar b)^{(0)},$
more precisely for some $C_0\subs C$ finite $\bar a$ is weakly isolated by $C_0\bar b$
over $C\bar b,$
 $\bar b\in \spn_Q(C_0\bar a),$ $C_0\bar b\le D$ and
$C_0$ weakly isolates $\bar b$ over $C.$ Let
    $\tau=(V,W,\{ W^{\frac{1}{l}}:l\in \N\})$ be a quantifier-free
type over $C_0$ satisfied by $\bar a\bar b,$ such that
$$\tau(\bar b)=(V(\bar b),W(\ex(\bar b)),\{ W^{\frac{1}{l}}(\ex(\bar b)):l\in \N\})$$
is a q-atom over $C\bar b$ and $\tau$ implies the q-atomic
type $\tau_0$ of $\bar b$ over $C.$
Then
$\bar x$ realises $\tau(\bar b)$ implies $\bar x\in \acl(C_0\bar b).$ Since
$\bar b\in \spn_Q(C_0\bar a),$ there is a $\Q$-linear combination with
parameters in $C_0$ \ $\bar q\cdot \bar x,$ such that
$\bar b=\bar q\cdot \bar a.$ Then the type saying that
$\bar x\bar y$ satisfies $\tau$ and $\bar y=\bar q\cdot \bar x$
is a q-atom of $\bar a$ over $C.$ Indeed, if $\bar a'$ satisfies the type
and $\bar b'=\bar q\cdot \bar a',$  then $\bar a'\bar b'$ satisfies $\tau.$
Then $\bar b'$ realises $\tau_0,$ thus $\bar b'$ is generic in $\tau_0$
over $C.$ It follows that $C\bar b'$ is isomorphic to $C\bar b.$ Then
$\tau(\bar b')$ is a q-atom over $C\bar b'$ satisfied by $\bar a'.$ It
implies that $\bar a'$ is generic in $\tau(\bar b')$ over $C\bar b'$, hence
$\bar a'\bar b'$ is generic in $\tau$ over $C.$
\qed
\bl If $C\le D,$ then\\
 For any $X\subs C^{(2n+2)}$
$$\trd(\ex(X)/\ex(C^{(2n+1)}))=\mdim(X/ C^{(2n+1)});$$
 For any $X\subs C^{(2n+1)}$
$$\trd(X/C^{(2n)})=\mdim(X/ C^{(2n)}).$$\el
\pf By definition of $\delta,$ using Lemma~\ref{ecl}(ii).\qed

\bl \lb{tilde0}  Assume $C\le A\le \D_1$ and $C\le B\le \D_2$  are
isomorphic over $C$ as substructures
by an isomorphism $\varphi:A\to B,$ $A=C+\spn(|\bar a|$ for a finite
string $\bar a,$
$C=\ecl(C)$ is countable and  the kernel in $\D_1,\D_2$ is standard.
  Then $\varphi$ can be
extended to an isomorphism $\varphi:\ecl(A)\to \ecl(B).$\el
\pf 
Claim 1.
If  $\varphi$ extends to $A\to B$ then it extends to
 $\acl(A)\to \acl(B).$
For this take any $\varphi$ preserving the field structure and notice that
for any $\Q$-linearly independent over $A$ elements
$a_1,\dots ,a_n\in \acl(A)$ the algebraic type of
$\ex(a_1),\dots ,\ex(a_n)$ over $\ex(A)$ is uniquely determined by the
requirement of algebraic independence.\\

Claim 2. $\varphi$  extends to $A+\spn(|\bar x|)$ for any $\bar x$ in
 $\ln(\acl(\ex(A))).$
Indeed, let $\ex(\bar x))$ be in $\acl(\ex(A)).$ Let
$W_A^{\frac{1}{l}}$ for $l\in \N$ be the
algebraic loci of $\ex({\frac{1}{l}}\cdot\bar v)$ over $\ex(A).$
$W_A$ can be represented as $W(\ex(\bar a))$ 
for some $W,$ an irreducible variety over $\ex(C),$ the latter algebraically
closed. We may assume $\bar x\bar a$ linearly independent over $C.$ 
By Lemma~\ref{1/l} the sequence $\{ W_A^{\frac{1}{l}}:l\in \N\}$
is determined by $W_A^{\frac{1}{l}}$ for some $l\in \N.$ Take
$\bar y$ in $\D$ such that
$$\ex({\frac{1}{l}}\cdot\bar y)\in W_B^{\frac{1}{l}}$$
where $W_B^{\frac{1}{l}}$ is obtained by applying
$\varphi$ to $W_A^{\frac{1}{l}}.$ Since $A,B\le D$ and
$$\trd(\ex(\bar x)/\ex(A))=0=\trd(\ex(\bar y)/\ex(B)),$$
necessarily $$\trd(\bar x/A)=n=\trd(\bar y/B),$$ i.e. both strings are
algebraically independent over the corresponding sets. Hence one can
extend $\varphi(\bar x)=\bar y.$\\
Claim 3.  If  $x\in \ecl(A),$ then there is $\bar x$ in 
$\ecl(A)$ extending $x,$
such that $A\bar x\le D_1$ and there is an isomorphism $A+\spn(|\bar x|)\to \ecl(B)$
extending $\varphi.$\\
Proof of Claim 3. Take for $\bar x$ a finite string from $\ecl(A)$ extending $x$
such that $A\bar x\le D.$
By Lemma~\ref{401}
$\bar x$ up to the ordering can be represented as $\bar x_0\dots\bar x_k$
such that $\bar x_{i+1}$ is in $\acl(A\bar x_0\dots\bar x_i)$ or
$\ex(\bar x_{i+1})$ is in $\acl(\ex(A\bar x_0\dots\bar x_i))$
depending on $i$ odd or even.\\

Since  $C$ is countable,  $\ecl(A)$ and
$\ecl(B)$ can be enumerated and by back-and-forth arguments using
Claim 3
one gets the isomorphism. \qed

\bc \lb{new_c} Suppose $\D_1, \D_2$ are with standard kernel, \ $A\le D_1$ and $B\le D_2$
are $\dd$-independent sets of cardinality $\lambda\le \omega_1,$  
$\varphi: A\to B$ is a bijection. Then $\varphi$ can be extended to an
isomorphism $\varphi:\ecl(A)\to \ecl(B).$\ec 
\pf 
First 
 enumerate $A$ and $B$  by ordinals agreeing with $\varphi.$ 
We can then define for each $i<\lambda$
$A_i\le A$  such that $$A_{i+1}= \{ a_j: j< i\}),$$ and correspondingly $B_i.$

Put
 $\varphi_0=\varphi_{|\ecl(\emptyset) }$ and
proceed by induction constructing an isomorphism  
$$\varphi_i: \ecl(A_i)\to \ecl(B_i).$$
Indeed, when a countable $\ecl(A_i)$ is extended by an element
we can apply
Lemma and back-and-forth arguments to construct $\varphi_{i+1}.$ 
Notice that, since $A$ is
$\dd$-independent, $\varphi\cup \varphi_i$ is still an isomorphism.

Take unions on limit steps.\qed

\bp \lb{tilde}  Assume $A\le \D_1$ and $B\le \D_2$  are isomorphic as substructures
by an isomorphism $\varphi:A\to B$ and  the kernel in $\D_1,\D_2$ is  compact.
  Then $\varphi$ can be
extended to an isomorphism $\varphi:\ecl(A)\to \ecl(B).$
\ep
\pf First we show that for any $x\in \ln(\acl(\ex(A)))$
$\varphi$ extends to $A+\spn(x)\to \D_2.$

Indeed, let
$W_A^{\frac{1}{l}}$ for $l\in \N$ be the
algebraic loci of $\ex({\frac{1}{l}}\cdot x)$ over $\ex(A).$
  By Lemma~\ref{n2} there is
$y\in D_2$ such that $$\ex({\frac{1}{l}}\cdot y)\in W_B^{\frac{1}{l}}$$
for all $l\in \N,$ where $W_B^{\frac{1}{l}}$ are obtained by applying
$\varphi$ to $W_A^{\frac{1}{l}}.$ Since $A\le \D_1,$ $B\le \D_2$ and
$$\trd(\ex(x)/\ex(A))=0=\trd(\ex( y)/\ex(B)),$$
necessarily $$\trd(x/A)=1=\trd(y/B),$$ i.e. both elements are
algebraically independent over the corresponding sets. Hence one can
define $$\varphi(a+{\frac{1}{l}}\cdot x)=b+{\frac{1}{l}}\cdot y$$ for
any $a\in A$ and $b=\varphi(a).$ \\
Now notice that Claim 1 of previous Lemma holds independently on the
kernel and the cardinality of $A.$ Combining it with the above proved
we can construct the required extension of $\varphi$ by transfinite induction 
and back-and-forth argument.\qed

\bl \lb{402}
Let $(V,W)$ be a free normal pair in $\D^n$  over some $C_0\le \D$
such that $\dim V +\dim W -n=d>0$ and $\dim V>1.$ Then for any $c_1,\dots,c_n\in D$
algebraically independent over $C_0$ there is $V'\subs V$ defined over
$C_0\bar c$ ($\bar c=\la c_1,\dots,c_n\ra$) such that $(V',W)$ is a normal
free pair and $\dim V' +\dim W -n=d-1.$ If $\bar a$ realises $(V',W)$ then
 $\delta(\bar c/C_0\bar a)<n.$\el
\pf Let $V'$ be an irreducible over $C_0\bar c$ component of the variety
$$V\cap \{ \bar x: c_1\cdot x_1 +\dots + c_n\cdot x_n=1\}.$$
Let $\bar a$ be a generic over $C_0\bar c$ point in $V'.$ Then $\bar a$
is generic in $V$ over $C_0$ and $\trd(\bar a/C_0\bar c)=\dim V-1,$
$\trd(\bar c/C_0\bar a)=n-1.$ In particular, it follows that
$\delta(\bar c/C_0\bar a)<n.$\\
{\bf Claim 1}. For any $k\le n$ and any distinct $i_1,\dots,i_k \in \{ 1,\dots,n\}$ \\
$\trd(a_{i_1},\dots,a_{i_k}/C_0)>\trd(a_{i_1},\dots,a_{i_k}/C_0\bar c)$
implies $\trd(a_{i_1},\dots,a_{i_k}/C_0)=\trd(\bar a/C_0).$\\
\pf Let $U\subs R^n$ be the minimal algebraic variety over
$C_0\cup \{ a_{i_1},\dots,a_{i_k}\} $ containing $\bar c.$ Suppose the
inequality holds. Then $\dim U<n.$ Also, $U$ must contain the hyperplane
$H_a$ defined by $\bar z\bar a=1,$ for otherwise $\trd(\bar c/C_0\bar a)<n-1.$
Hence $H_a$ is a component of $U.$ But then
$\bar a\in \acl(C_0\cup \{ a_{i_1},\dots,a_{i_k}\}),$ since $\bar a$ is the only
vector satisfying $\bar z\bar a=1$ for all $\bar z\in H_a.$ It follows that
$\trd(a_{i_1},\dots,a_{i_k}/C_0)=\trd(\bar a/C_0).$\\
{\bf Claim 2}. For any non-zero integer $n$-tuple $\bar m=\la m_1,\dots,m_n\ra$
the element $\bar m\cdot\bar a=m_1\cdot a_1+\dots+m_n\cdot a_n$ is not in
$\acl(C_0\bar c).$\\
\pf Notice first that $a\notin \acl(C_0),$ since $V$ is free
of additive dependencies. Thus, if $a\in \acl(C_0\bar c),$
then $\trd(\bar c/C_0a)<n.$ By the argument from the proof
of Claim 1 it would follow that $|\bar a|\subs \acl(C_0a).$
This contradicts the fact that $\dim V>1.$\\
\pf of Lemma. It follows from Claim 2, that $V'$ is free of additive dependencies.
Thus the pair $(V',W)$ is free. To prove the normality of the pair consider
for any distinct $i_1,\dots,i_k \in \{ 1,\dots,n\}$ the number
$$\dim V'_{i_1,\dots,i_k} +\dim W_{i_1,\dots,i_k} -k.$$
Suppose towards the contradiction that the number is negative. Then
$$\dim V_{i_1,\dots,i_k} +\dim W_{i_1,\dots,i_k} -k\le 0,$$
$\dim V'_{i_1,\dots,i_k}\ge \dim V_{i_1,\dots,i_k} -1.$ It follows that
$\dim V'_{i_1,\dots,i_k}< \dim V,$ since
$\dim W\le \dim W_{i_1,\dots,i_k}+(n -k)$ and $d>0.$ Thus
$$\trd(a_{i_1},\dots,a_{i_k}/C_0)<\trd(\bar a/C_0).$$
Then by Claim 1
$$\trd(a_{i_1},\dots,a_{i_k}/C_0)=\trd(a_{i_1},\dots,a_{i_k}/C_0\bar c).$$
Thus $\dim V'_{i_1,\dots,i_k}= \dim V_{i_1,\dots,i_k}$ and hence
$$\dim V_{i_1,\dots,i_k} +\dim W_{i_1,\dots,i_k} -k=
\dim V'_{i_1,\dots,i_k} +\dim W_{i_1,\dots,i_k} -k<0,$$
contradicting the normality of $(V,W).$

To finish the proof of normality we need to consider
also
$$a'_i=m_{i,1}a_1+\dots+ m_{i,n}a_n,\ \ \mbox{ and }
b'_i=b_1^{m_{i,1}}\cdot\dots \cdot b_n^{m_{i,n}}\ \ i=1,\dots n$$
where $\la m_{i,1},\dots, m_{i,n}\ra$ are independent integer vectors
and
$\la b_1,\dots, b_n\ra$ a generic element of $W$ over 
$\ex(C_0\bar c).$ 
The corresponding inequalities for the loci over
$C_0\bar c$ and $\ex(C_0\bar c)$ of the $\bar a'$ and $\bar b'$ are
proved by the same argument.   
\qed
\bl \lb{L1} Let $C\le D$   and $V,V'\subs D^n$   algebraic varieties
  over a finite $C_0\subs C,$ \ $W, W'\subs R^n$ varieties  over $\ex(C_0).$
Let $\bar a$ in $\D$ realise the type  $\tau=(V\smin V',W\smin W',\{ W^{\frac{1}{l}}:l\in \N\})$
and satisfy the following minimality condition:\\
for any $\bar a'\in \D$ realising $\tau$
either $\mdim(\bar a'/C)>\mdim(\bar a/C)$
or $\mdim(\bar a'/C)=\mdim(\bar a/C)$ and
$\delta(\bar a'/C)\ge \delta(\bar a/C).$

Suppose the ex-locus   $(V_0,W_0,\{ W_0^{\frac{1}{l}}:l\in \N\})$
of $\bar a$ over $C$ is weakly defined over $C_0.$ Then 
$\tau_0=(V_0\smin V',W_0\smin W',\{ W_0^{\frac{1}{l}}:l\in \N\})$
is a q-atom over $C$ weakly defined over $C_0.$ If $C=\ecl(C),$ then
$\delta(\bar a/C)=0$ and there is a finite $C_1\subs C$ such that $\tau$ is
defined over $C_1.$\el
\pf Let $\bar a'$  be any realisation of $\tau_0$ and assume
towards a contradiction that $\bar a'$ is not generic over $C.$ Then
$$\trd(\bar a'/C)+\trd(\ex(\bar a')/\ex(C))<
\trd(\bar a/C)+\trd(\ex(\bar a)/\ex(C)).$$
 By definition $\mdim(\bar a'/C)\ge \mdim(\bar a/C),$
so  $\delta(\bar a'/C)<\delta(\bar a/C),$ the contradiction.
This proves the first part of the statement.

Assume $C=\ecl(C)$ and assume towards a contradiction that
$\delta(\bar a/C)=d>0.$
We can reduce the problem to a $\Q$-linear basis of
$\bar a$ over $C,$ so
we  assume that $\bar a=\la a_1,\dots,a_n\ra$ is  $\Q$-linearly independent
over $C.$ It follows that  $(V_0,W_0)$ is absolutely free and normal.

By the assumptions $\dim V_0+\dim W_0-n=\delta(\bar a/C)>0.$
Take $c_1,\dots,c_n\in C$ algebraically independent over $C_0$ (such
ones exist
in $C_0^{(1)}).$
By Lemma~\ref{402} either $\dim V_0=1,$ or there is $V'_0\subs V_0$ over
$\bar C_0\bar c$ such that $(V'_0,W_0)$ is free and normal. In the second case
by Lemma~\ref{ec1} there is a generic over $C_0\bar c$ realisation
$\bar a'$ of
$$(V'_0,W_0,\{ W_0^{\frac{1}{l}}:l\in \N\}),$$ which is also generic
realisation of $\tau_0$ over $C_0,$ such that
$$\delta(\bar a'/C_0\bar c)=\dim V'_0+\dim W_0-n<d.$$
The contradiction.
In the first case, since $d>0,$ necessarily $\dim W_0=n,$ i.e. $W_0=R^n.$
Choose $c\in C\smin C_0^{(1)}$ and let
 $$W'_0=\{ \la y_1,\dots,y_n\ra: y_1+\dots+ y_n=c\}.$$
Since $\ex(c)\notin \acl(\ex(C_0)),$ any generic over $\ex(C_0c)$ realisation
$\la b_1,\dots,b_n\ra$ of $W'$ is generic in $W$ over $\ex(C_0).$
Also, if $n>1,$
for any $k<n$ and distinct $i_1,\dots,i_k \in \{ 1,\dots,n\}$
$$\trd(b_{i_1},\dots,b_{i_k}/\ex(C_0c)) =\trd(b_{i_1},\dots,b_{i_k}/\ex(C_0)) ,$$
which yields the normality of $(V_0,W'_0).$  The same argument shows
that\\ $b_1^ {m_1}\cdot\dots \cdot b_n^{m_n}\notin \acl(\ex(C_0c))$ for any
non-zero integer $n$-tuple $\la {m_1},\dots,{m_n}\ra,$ which shows that
$W'_0$ is free of multiplicative dependencies. Thus $(V_0,W'_0)$ is free.
Again it gives by Lemma~\ref{ec1} as above a generic over $C_0$
realisation of $\tau_0$ over $C_0,$ such that
$$\delta(\bar a'/C_0\bar c)=\dim V_0+\dim W'_0-n<d.$$ The final
contradiction.

Finally the associated sequence is defined over some finite $C_1$ by Corollary~\ref{c2}.

\qed

\bl \lb{ec2} Let  $\ecl(\emptyset)\subs C\le D,$
$C\smin \ecl(\emptyset)$ be finite,  $\bar a\bar b$ a tuple in
$\D'\ge \D,$ $(V\smin V_1,W\smin W_1,\{ W^{\frac{1}{l}}:l\in \N\})$ be 
be a type over $C$ realised by
$\bar a\bar b,$ $(V_0,W_0,\{ W_0^{\frac{1}{l}}:l\in \N\}),$ $(V,W)$
the main part of ex-locus of the tuple,  and assume $\delta(\bar a/C)=0.$
Then there is $(V_0',W_0')$ over $C$ such that for any $\bar a'$
realising 
$\sigma=(V_0\smin V_0',W_0\smin W_0',\{ W_0^{\frac{1}{l}}:l\in \N\})$ \
there is $\bar b'$ in some $\D''\ge \D'$ such that $\bar a'\bar b'$
realise $\tau.$ \el
\pf  
We just let $V_0\smin V_0'$ be the projection of $V$ onto the
coordinate space corresponding to $\bar a$ and $W_0\smin W_0'$ 
the projection of $W.$

We prove the statement by induction on $|\bar b|=n\ge 1$ assuming also
that $\bar a\bar b$ are $\Q$-linearly independent over $C.$ It follows from
the assumptions that $(V,W)$ is an absolutely free and normal pair.

Let $n=1.$ If both
 $\dim V(\bar a)=1$ and $\dim W(\ex(\bar a))=1,$ then, by the elementary
algebraic geometry, for any
 $\bar a'$ realising $\sigma$
 $\dim V(\bar a')=1$ and $\dim W(\ex(\bar a'))=1.$ Hence $V(\bar a')$ and
$W(\ex(\bar a'))$ are just the whole affine lines and we can take
for $b'$ any point in $D.$\\
If  $\dim V(\bar a)=0,$ then necessarily $\dim W(\ex(\bar a))=1$ and
again $\dim V(\bar a')\ge 0$ and $\dim W(\ex(\bar a'))=1.$ Take
for $b'$ any point in $V(\bar a').$\\
If $\dim W(\ex(\bar a))=0,$ then $\dim V(\bar a)=1$ and
$\dim W(\ex(\bar a'))\ge 0,$  $\dim V(\bar a')=1.$

Let $$\tau_a=(V(\bar a\bar c),W(\ex(\bar a\bar c)),\{ (W(\ex(\bar a\bar c)))^{\frac{1}{l}}:l\in \N\})$$ be the
quantifier-free type of $\bar b$ over $C\bar a,$ where $\bar c$ is a finite string
from $C$ and $V,$ $W$ are irreducible varieties over $\ecl(\emptyset),$
$\ex(\ecl(\emptyset)),$ correspondingly.

 $\tau_{a'}$ is the type obtained by replacing all occurencies
of $\bar a$ by $\bar a'.$ By Lemma~\ref{1/l} the sequence
$\{ (W(\ex(\bar a'\bar c)))^{\frac{1}{l}}:l\in \N\})$
 is determined by its cut of a length
$l_0,$ thus there is $y$ in $R$ and an associated sequence such that
$$y^{\frac{1}{l}}\in
(W(\ex(\bar a'\bar c)))^{\frac{1}{l}}\mbox{ for all }l\in \N\}.$$
Take $b'\in D$ such that
$$\ex({\frac{1}{l_0}}\cdot b')=y^{\frac{1}{l_0}}.$$
Then
$$\ex({\frac{1}{l}}\cdot b')\in (W(\ex(\bar a'\bar c)))^{\frac{1}{l}}\mbox{ for all }l\in \N.$$
Since $V(\bar a')=D,$ \ $b'$ satisfies $\tau_{a'}$ and $\bar a'b'$ satisfies $\tau.$

For $n>1$  notice first that by the assumptions  $C\bar a\le D,$ thus
$(V(\bar a),W(\ex(\bar a)))$ is a normal pair. To prove the statement of
Lemma we consider two alternative cases:\\
Case 1. There are distinct ${i_1},\dots {i_k}$ in
$\{ 1,\dots,n-1\}$ such that the projections
$V_{{i_1},\dots {i_k}} (\bar a) $ and $W_{{i_1},\dots {i_k}}(\ex(\bar a))$
onto corresponding coordinates satisfy
$$\dim V_{{i_1},\dots {i_k}} (\bar a) +
\dim W_{{i_1},\dots {i_k}}(\ex(\bar a)) =k.$$
Let for simplicity $\la {i_1},\dots {i_k}\ra =\la 1,\dots k\ra.$
Then $\delta(\la b_{1},\dots b_{k}\ra/C\bar a)=0$ and  one can consider
$\bar a_+=\bar a\la b_{1},\dots b_{k}\ra$ in place of $\bar a$ and
$\bar b_-=\la b_{k+1},\dots b_{n}\ra$ in place of $\bar b.$
By induction hypothesis  find type $\sigma_+$ over $C$ which says
about any $\bar a'_+$ that there is $\bar b'_-$ such that
$\bar a'_+\bar b'_-$ satisfy $\tau.$  \\
Case 2.  $$\dim V_{{i_1},\dots {i_k}}(\bar a)+\dim W_{{i_1},\dots {i_k}}(\ex(\bar a))>k$$
for any distinct ${i_1},\dots {i_k}\in \{ 1,\dots,n-1\}.$ \\
If the pair $(V(\bar a), W(\ex(\bar a))$ is free then the statement follows
from Lemma~\ref{ec1}.\\
Suppose  $W(\ex(\bar a))$ has a multiplicative dependence:
$$y_1^{m_1}\cdot\dots\cdot y_n^{m_n}=r\in \acl(C\bar a)$$
   for some ${m_1},\dots {m_n}\in\Z,$
$r\in R.$ Notice that by definitions
 $$\ex(b_1)^{m_1}\cdot\dots\cdot \ex(b_n)^{m_n}=r$$ holds, and so
$$m_1{b_1}+\dots\cdot m_n{b_n}=d$$ holds for some $d\in \D$ such that
$\ex(d)=r.$   Assume $m_n\neq 0.$ Apply the induction hypothesis
considering $\bar ad$ in place of $\bar a$ and $\la b_1,\dots,b_{n-1}\ra$
in place of $\bar b.$ We have then
that whenever $\bar a'd'$ satisfy the same type as $\bar ad,$ there is
$\la b'_1,\dots,b'_{n-1}\ra$ such that $\bar a'\la d',b'_1,\dots,b'_{n-1}\ra$
is of the same type as
$\bar a\la d,b_1,\dots,b_{n-1}\ra.$ On the other hand, considering $n=1,$ we proved that
whenever $\bar a'$ satisfies the type of $\bar a,$ there is $d'$ such that
$\bar a'd'$ satisfy the type of $\bar ad,$ which completes the proof in the
case.

If  $V(\bar a)$ has an additive dependence we act symmetrically. \qed
\bl \lb{L2} Let $\ecl(\emptyset)\subs C\le D$,  \ $\bar a,$ $\bar b$ in $\D,$
$\delta(\bar a/C)=0$ and $\bar a$ q-atomic
over $C,$ $\bar b$ q-atomic over $C\bar a.$ Then $\bar a\bar b$ is q-atomic
over $C.$\el
\pf Choose finite $C_0\subs C$ such that $C_0$ weakly isolates $\bar a$ over $C$ and
$C_0\bar a$ weakly isolates $\bar b$ over $C\bar a.$
Let $\tau=(V\smin V_1,W\smin W_1, \{ W^{\frac{1}{l}}:l\in \N\})$ be
the q-atom of 
$\bar a\bar b$ over $C_0$ and
$\sigma$ the type obtained by projecting  the
varieties  onto $\bar a$-coordinates. Then
$\sigma$ is a q-atom weakly isolating $\bar a$ over $C.$
Also, the type of $\bar b$ over  $C\bar a$ is determined over
$C_0\bar a$ by $\tau(\bar a).$
We claim that whenever $\bar a'\bar b'$ in some $\D'\ge \D$
 realises $\tau,$   \ $\bar a'\bar b'$ is
generic over $C.$ This would prove Lemma.

Assume  $\bar a'\bar b'$ realises $\tau.$
Then  $\bar  a'$  satisfies $\sigma$ and, since $\sigma$
is a q-atom over $C,$ the quantifier-free types of $\bar a$ and $\bar a'$
over $C$ coincide.

Suppose towards a contradiction that $\bar b'$ is not generic in $\tau(\bar a')$
over $C.$ Then $\bar b'$ satisfies $\tau(\bar a')$ and a pair
$(V'(\bar a'),W'(\ex(\bar a')))$
with $$\dim V'(\bar a')< \dim V(\bar a')\mbox{ or }
\dim W'(\ex(\bar a'))< \dim W(\ex(\bar a')).$$
Since $\tau$ and $(V',W')$ are defined over some finite
set $C_1,$ $C_0\subs C_1\subs C,$ by Lemma~\ref{ec2}  there is $\bar b''$ which satisfies
$\tau(\bar a)$
and $(V'(\bar a),W'(\ex(\bar a)).$ This contradicts the assumption that
$\tau(\bar a)$ is a q-atom over $C.$

  \bl \lb{L3} Let $\ecl(\emptyset)\subs C\le C\bar a\le D$ and $\bar a\bar b$ be q-atomic over $C.$
 Then $\bar a$ is q-atomic over $C$ and $\bar b$ is q-atomic over $C\bar a.$\el
\pf  Let
$$\tau=(V\smin V',W\smin W',\{ W^{\frac{1}{l}}:l\in \N\})$$  be the q-atomic type
of $\bar a\bar b$ over $C$ weakly defined over some finite $C_0.$
Let $V_1\smin V_1',W_1\smin W_1'$ be the projections of $V\smin
V',W\smin W'$ onto the coordinates corresponding
to $\bar a,\ex(\bar a)$ and
$$\sigma=(V_1\smin V_1',W_1\smin W_1',\{ W_1^{\frac{1}{l}}:l\in \N\}).$$
Suppose $\bar a'$ realises type $\sigma.$ Then,
by Lemma~\ref{ec2}, there is $\bar b'$ such that  $\bar a'\bar b'$ realises
$\tau.$ Since $\tau$ is a q-atom over $C,$
$\bar a'\bar b'$ is generic in $(V,W)$ over $C.$ It follows $\bar a'$
is generic in $(V_1,W_1)$ and $\bar b'$ generic in
$(V(\bar a'),W(\ex(\bar a'))).$ The first fact yields that $\bar a$ is
q-atomic over $C,$ the second fact, after assuming $\bar a'=\bar a,$ yields
$\bar b$ is q-atomic over $C\bar a.$\qed
\df A structure $A\in \SED$  is said to be {\bf q-atomic over $C\le A$}
 if for any
$\bar a$ in $A$ such that $C\bar a\le A$ the tuple $\bar a$ is q-atomic over
$C.$\\ \\

\bl \lb{L5} Assume $A$ is with full kernel, $\ecl(\emptyset)\subs C\le A.$

(i) If $B\subs A$ finite, $\delta(B/C)=0$ and $A$ is q-atomic over $C,$ then $A$ is q-atomic
over $CB.$

(ii) If $\ecl(C)\subs A$ and $A$ is q-atomic over $\ecl(C),$ then
$A$ is q-atomic over $C$ and for any $\bar a$ in $A,$ such that $C\bar a\le A,$
it holds $\delta(\bar a/C)=0.$
\el
\pf (i) Suppose $CB\bar a\le A.$ By assumtions, $B\bar a$ is q-atomic over $C.$
Then by Lemma~\ref{L3} $\bar a$ is q-atomic over $CB.$ \\
(ii) Suppose $C\bar a\le A.$ Extend $\bar a$ to $\bar a'$ so, that
$\ecl(C)\bar a'\le A.$ Then, by assumtions, $\bar a'$ is q-atomic over $\ecl(C)$
and by \ref{L1} $\delta(\bar a'/\ecl(C))=0.$ \\
Choose finite $\bar c$ in $\ecl(C)$ such that $\bar a'$ is q-atomic over
$C\bar c,$ $\delta(\bar a'/C\bar c)=0$ and $C\bar c\le A.$ By Lemmas~\ref{401}
and \ref{L2}, $\bar c\bar a'$ is q-atomic over $C.$ Also
$$\delta(\bar c\bar a'/C)=\delta(\bar a'/C\bar c)+\delta(\bar c/C)=0.$$
Since $\bar c\bar a'$ extends $\bar a$ and $C\bar a\le A,$ \
$\delta(\bar a/C)=0.$ Finally, by Lemma~\ref{L3}, $\bar a$ is q-atomic over $C.$
\qed

\bp \lb{L4} Over
any $C\le \D$   there is a q-atomic structure $E(C)\in \EC,$ $E(C)\le D.$ 
If the kernel $K$ is compact then
 for any embedding $C\le \D'\in \EC$  there is an extension of the
embedding $E(C)\le \D'.$
  
The cardinality of $E(C)$ is not greater than
$\card(C)+\card(K).$
\ep
\pf  
$E_D(C)=E(C)$ will be represented as  $\bigcup_n C_n$ for some
$$C_0\le \dots \le C_n\le \dots\subs \D.$$
Put $C_0=\ecl(C)$ and assume $C_n$ is constructed. Let
$\{ \mu_{n,\alpha}: \alpha<\lambda_n\}$ be the set of all
types of the form $(V, W),$ irreducible free and normal
over $C_n,$  $\lambda_n=\card(C_n).$

Put $C_{n,0}=C_n$ and, when $C_{n,\alpha}$ is constructed,
use Lemma~\ref{L1} to find in some normal extension of $C_{n,\alpha}$
a q-atomic over $C_{n,\alpha}$
tuple $\bar a_{n,\alpha}$ realising $\mu_{n,\alpha}.$ Since, by
Lemma~\ref{normal} the incomplete ex-locus of $\bar a_{n,\alpha}$ over $\ecl_D(C_{n,\alpha})$ is
absolutely free and normal, by Lemma~\ref{almost} the corresponding type is
realisable in  $\D,$ so we assume the tuple is in
$\D.$  

Put $$C_{n,\alpha+1}=C_{n,\alpha}\cup |\bar a_{n,\alpha}|.$$
For limit ordinals $\sigma$
$$C_{n,\sigma}=\bigcup_{\alpha <\sigma}C_{n,\alpha}.$$
Finally, $C_{n+1}=C_{n,\lambda_n}.$\\
{\bf Claim}. For any $n$ and $\alpha<\lambda_n$

(i)  $C_{n,\alpha}$ is q-atomic over $C_0;$

(ii) $C_{n,\alpha}\le D;$

(iii)  $C_{n,\alpha}$ is q-atomic over $C;$

(iv) if kernel is compact,  $C_{n,\alpha}$ does not depend on $\D.$

\pf of Claim. By induction on lexicographically ordered pairs
$\la {n,\alpha}\ra. $ \\
For $\la 0,0\ra$ the statements follow from properties of $\ecl.$ The inductive
step for limit ordinals is trivial, so we assume the Claim is true for
$C_{n,\alpha}$ and want to prove it for $C_{n,\alpha+1}.$\\
Let $\bar a$ be a tuple from $C_{n,\alpha+1},$ such that
$C_0\bar a\le D.$
By the construction there is $\bar b$ in $C_{n,\alpha}$ which weakly isolates
$\bar a_{n,\alpha}$ over $C_{n,\alpha}.$ When we extend $\bar b$ the property
is preserved, so we assume $C_0\bar b\le D$ and
$|\bar a|\subs |\bar b\bar a_{n,\alpha}|.$ Then, by induction, $\bar b$
is q-atomic over $C_0$ and by Lemma~\ref{L1} $\delta(\bar b/C_0)=0.$
Now by Lemma~\ref{L2} $\bar b\bar a_{n,\alpha}$ is q-atomic over $C_0.$
Since $|\bar b\bar a_{n,\alpha}|=|\bar a\bar b'|$ for some $\bar b',$
by Lemma~\ref{L3} $\bar a$ is q-atomic over $C_0,$ which proves (i).

To prove (ii) we must show that $\delta(X/C_{n,\alpha+1})\ge 0$ for any
finite $X\subs D.$ Notice that
$\delta(X/C_{n,\alpha+1})=\delta(X/C_0Y)$ for large enough finite
$Y\subs C_{n,\alpha+1}.$ We can choose $Y=|\bar b\bar a_{n,\alpha}|$ with
$\bar b$ as above. Then $Y$ is q-atomic over $C_0$ and by Lemma~\ref{L1}
$\delta(Y/C_0)=0.$ Hence $\delta(X/C_0Y)\ge 0,$ by Lemma~\ref{tr}.
Which finishes the proof of (ii).

(iii) follows from (i) and Lemma~\ref{L5}.

(iv) follows from Proposition~\ref{tilde}.
 
It follows from the (i)-(iii) of the claim that $E(C)$ is q-atomic over $C.$ It follows
from Theorem\ref{t1} that
$E(C)$ is e.a.c.  In case the kernel is compact it 
follows from (iv) of the claim that $E(C)$ can be embedded in any
$\D'$ provided $C\le \D'.$

The cardinality statement follows from the fact that
$\lambda_n\le \card(C_0)= \card(\ecl(C)\le \card(C)+\card(K).$ 
\qed

\df A structure $E(C)\in \EC$  is said to be
{\bf constructible q-prime over $C$} if

(i) $C\le E(C),$

(ii) $E(C)$ is q-atomic over $\ecl(C),$

(iii) there is a sequence $\{ \bar a_j:j< \lambda \}$ of tuples of $E(C)$
for some limit ordinal $\lambda,$ such that
$E(C)=\bigcup\{ |\bar a_j|:j< \lambda \}$ and
for any $j\le\lambda$   \ \
$\bar a_j$ is
a realisation of a q-atom $\tau_j$ over $\ecl(C)\cup A_j,$ where
$A_{j}=\bigcup_{i<j}|\bar a_{i}|.$ The corresponding sequence
$\{ (\bar a_j,\tau_j):j< \lambda \}$  is said to be a {\bf q-construction}.
We assume $|\bar a_j|\cap(\ecl(C)\cup A_j)=\emptyset.$\\ \\

\bp \lb{p2} Suppose the kernel $K$ is compact. Then over
any $C\in \SED$   there is a q-prime structure $E(C)$ such
that for any embedding $C\le \D\in \EC$  there is an extension of the
embedding $E(C)\le \D.$ 
Any two constructible q-prime structures over $C$  are isomorphic over  $C.$\ep
\pf The first statement is proved in Proposition~\ref{L4}, where the sequence
$\{ \bar a_j:j<\lambda\} $ is just
$\{ \bar a_{n,\alpha}:n<\omega,\alpha<\lambda_n\} $ getting ordinal
enumeration and $\tau_j$ are the corresponding types.
  To prove the rest we use the argument
from [Sh]
following [B].\\

We assume by \ref{tilde} that $C=\ecl(C).$

A subset $B$ of a constructible q-prime structure over $C$
is said to be closed if for all $\bar a_j,$ if
$|\bar a_j|\cap B\neq \emptyset,$ then $|\bar a_j|\subs CB,$ \
$CB$ contains all parameters of type $\tau_j$ and $CB\le E(C).$
Notice that each point $b\in B$ can
be included in some finite string $\bar b$ from $B$ such that
$C\bar b\le E(C).$ Since $E(C)$ is q-atomic, by Lemma~\ref{L1}
$\delta(\bar b/C)=0.$ In particular, if $B$ is finite $\delta(B/C)=0.$ \\
{\bf Claim 1}. If $B$ is a closed subset of $E(C),$ $j<\lambda,$ then

(i) $CBA_j\le E(C);$

(ii) $\delta(\bar a_j/CBA_j)=0;$

(iii) $\bar a_j$ is q-atomic over $CBA_j$ and the corresponding q-atom
is $\tau_j.$\\
Proof of Claim.  For any $j<\lambda$  we now prove (i) - (iii)
for $B^i=B\cap A_i,$  in place of $B$ for all $i\le \lambda$ by induction
on $i.$ \\
For $i\le j$ this is true, since then $B^i\subs A_j,$ (i) and (iii) are true
 by the
definition and (ii) follows by the following argument, which works for any
$i,$ provided (i) and (iii) hold. Choose finite
$A\subs B^iA_j,$ such that $CA\le E(C)$ and $A$ contains the parameters of
$\tau_j.$ Then $A$ is q-atomic over $C$ and $\bar a_j$ is q-atomic over
$CA.$ By Lemma~\ref{L3} $A\bar a_j$ is q-atomic over $C,$ hence by 
Lemma~\ref{L1} \ $\delta(A\bar a_j/C)=0,$ $\delta(A/C)=0.$ It follows
$\delta(\bar a_j/CA)=0,$ which in turn implies $\delta(\bar a_j/CB^iA_j)=0.$\\
Assume now that (i)-(iii) hold for some $i,$ $j\le i<\lambda,$ we want to
prove it for $i+1.$ Assume $B^{i+1}\neq B^i.$ Then
$B^{i+1}=|\bar a_{i}|\cup B^{i}.$ Since $B^iA_j\subs A_i$ and the parameters of
$\tau_i$ are in $B^i,$ \ $\bar a_i$ is q-atomic over $CB^iA_j\bar a_j.$
By the induction hypothesis $\bar a_j$ is q-atomic over $CB^iA_j$ by $\tau_j$
and
$\delta(\bar a_j/CB^iA_j)=0.$ Then by \ref{L2} $\bar a_j\bar a_i$ is q-atomic
over $CB^iA_j.$ More precisely, the type of $\bar a_j\bar a_i$ over
$CB^iA_j$ is determined by $\tau_j(\bar x)\&\tau_i(\bar y)$ (where $\bar x,$
$\bar y$ correspond to $\bar a_j,$ $\bar a_i$).\\
Since $\bar a_i$ is q-atomic over $CB^iA_j$ and $CB^iA_j\le E(C),$ by the
argument above $\delta(\bar a_i/CB^iA_j)=0.$ Then $CB^iA_j\bar a_i\le E(C).$
Hence, by Lemma~\ref{L3} and the above proved,
 $\bar a_i$ is q-atomic  over $CB^iA_j$ and $\bar a_j$ is q-atomic over
$CB^iA_j\bar a_i,$ more presicely the q-atom is again $\tau_j.$
Since $CB^iA_j\bar a_i=CB^{i+1}A_j,$ we proved (i) and (iii)
for $B^{i+1}.$ Once again (ii) is proved by the argument above. \\
For limit ordinals $i$ (i) - (iii) is evident. This finishes the proof of
Claim.\\
{\bf Claim 2}. Suppose $B$ is closed. Then for any finite $X\subs E(C)$
there is finite $B_X\subs E(C),$ such that $B B_X$ is closed and
$X\subs B_X.$\\
Proof of Claim. We prove by induction on $j$ that for
$X\subs A_{j}$ there is $B_X\subs A_j,$ such that $BB_X$ is closed.
If $j=0,$ $X=\emptyset=B_X.$ Suppose  $|\bar a_j|\subs X\subs A_{j+1}.$
Choose finite $Y\subs A_j,$ such that $X\cap BA_j\subs Y,$
$CBB_Y$ is closed and $CB_Y$ contains all the parameters of $\tau_j.$
Then there is finite $Z\subs B$ such that $CZB_Y\le CBB_Y\le E(C)$
and $\bar a_j$ is q-atomic over $CZB_Y.$ Since $E(C)$ is q-atomic over $C,$
$ZB_Y$ is q-atomic and by \ref{L1}, \ref{L2} $ZB_y\bar a_j$ is q-atomic over
$C.$ Hence $\delta(ZB_y\bar a_j/C)=0,$ thus $CZB_Y\bar a_j\le E(C).$
Since $Z$ can be chosen as large as we want, it follows $CBB_Y\le E(C).$
We state that $BB_Y\bar a_j$ is closed. It remains to prove that
if there is $a\in |\bar a_i|$ such that $a\in BB_Y|\bar a_j|$ then
$|\bar a_i|\subs CBB_Y|\bar a_j|$ and the parameters of $\tau_i$ are
in the same set. By the construction it is true if $a\in BB_Y.$ If
$a\in |\bar a_j|,$ then $i=j$ and both $|\bar a_j$ and the set of all
parameters of $\tau_j$ are subsets of $CBB_Y|\bar a_j|.$ Thus
$BB_Y|\bar a_j|$ is a closed set which contains $X.$ Claim proved.\\
{\bf Claim 3}. $E(C)$ is q-atomic over $CB$ for any closed $B$
and if $CB\bar a\le E(C),$ then $\delta(\bar a/CB)=0.$\\
Proof of Claim. It is enough to prove that any $\bar a$ in $CBA_j,$ such that
$CB\bar a\le E(C),$ \ $E(C)$ is q-atomic over $CB$ and $\delta(\bar a/CB)=0$
$CB$ for all $j<\lambda.$ We prove it by induction on $j.$ The statement
is trivial for $j=0.$ Suppose it holds for $j$ and $\bar a$ is a tuple in
$CBA_{j+1}$ such that $CB\bar a\le E(C).$ Choose finite tuple $\bar b $
in $CBA_j$ such that $CB\bar b\le E(C),$ $|\bar a|\subs |\bar b\bar a_j|$
and $\bar a_j$ is q-atomic over $C\bar b$ (Claim 1). By Lemma~\ref{L2}
$\bar b\bar a_j$ is q-atomic over $CB.$ By Lemma~\ref{L3} $\bar a$ is q-atomic
over $CB.$ Since      $|\bar b\bar a_j|=|\bar a\bar b'|$ for some $\bar b',$
by Lemma~\ref{L3} $\bar a$ is q-atomic over $CB.$
For  large enough $\bar b$
$$\delta(\bar b\bar a_j/CB)=\delta(\bar a_j/CB\bar b)+\delta(\bar b/CB)=
\delta(\bar a_j/CB\bar b)=\delta(\bar a_j/CBA_j)=0.$$
Since $$0=\delta(\bar a\bar b'/CB)=\delta(\bar b'/CB\bar a)+\delta(\bar a/CB)$$
and both summands on the right are non-negative by $CB\le CB\bar a\le E(C),$
we get $\delta(\bar a/CB)=0.$ Claim proved.\\

Now we come to the proof of Proposition. We consider two constructible
q-prime structures $A$ and $A'$ over $C$ with q-constructions
$\{ (\bar a_j,\tau_j):j< \lambda\}$ and $\{ (\bar a'_j,\tau'_j):j< \lambda'\}$
 correspondingly. Assume $\lambda \le \lambda'.$ We define by induction
ascending chains of isomorphic closed sets  $ C_j\le A$ and $ C'_j\le A'$
for $j<\lambda$ and the isomorphisms  $f_j:C_j\to C'_j.$ We also satisfy the
condition that if  $j=\eta+(2n+2)$    $\eta$ limit, $n$ natural, then
$\bar a_{\eta+n}$ lives in $C_j;$ if $j=\eta+(2n+1)$ then
$\bar a'_{\eta+n}$ lives in $C'_j.$\\
Consider a typical case of the construction: $j=\eta+(2n+2).$ By
Claim 2 there is a finite set $B_0$ containing $|\bar a_{\eta+n}|,$
such that $C_{j-1}B_0$ is closed. By Claim 3 $A$ is q-atomic over $C_{j-1}$
and $\delta(B_0/C_{j-1})=0.$ Since
$A$ is q-atomic over $C_{j-1}$ there is an isomorphism $g_0$ extending
$f_{j-1}$ and taking $B_0$ onto some $B'_0\subs A'.$
By Claim 2
$B'_0$ can be included
in some  $B'_1\subs A'$ such that $C'_{j-1}\cup B'_1$ is closed.
By the isomorphism $\delta(B'_0/C'_{j-1})=0,$ hence by Lemma~\ref{L5}
  $A'$ is q-atomic over
$C'_{j-1}B'_0.$ Thus $B'_1$ satisfies a q-atom over $C'_{j-1}B'_0.$
It follows there is an isomorphism $g_1$ extending
$g_0$ and taking $C_{j-1}\cup B_1$ for some  $B_1$ onto $C'_{j-1}\cup B'_1.$
Continuing in this way we get a chain of finite sets $B_m\subs A$ ($m\in \N$),
  $B'_m\subs A'$ and corresponding  isomorphisms $g_m$ from
$C_{j-1}B_m$ onto $C'_{j-1}B'_m,$ such that $C_{j-1}B_m$ is
closed if
$m$ is even and $C'_{j-1}B'_m$ is closed if $m$ is odd. Finally the sets
$$C_j=C_{j-1}\cup \bigcup_{m\in \N}B_m,  \ \ \
C'_j=C'_{j-1}\cup \bigcup_{m\in \N}B'_m\mbox{ and }f_j=\bigcup_{m\in \N} g_m$$
satisfy the conditions required for $j.$  For odd $j$ the procedure is
symmetrically from $A'$ to $A.$   In the end the construction yields an
isomorphism $f=f_{\lambda}$ from $A$ onto the closed  subset $C'_{\lambda}$
of $A'.$ Since $A'$ is q-atomic over $C'_{\lambda}$ and by the isomorphism
any q-atom over $C'_{\lambda}$ is  realised in $C'_{\lambda}$ we have
necessarily $C'_{\lambda}=A'.$\qed
\bt Let $C\le D$ and $C$ be a $\dd$-independent set of power
$\lambda$ and $\kappa=\card K.$ Suppose that $D$ is  q-prime constructible over $C.$ Then

(i) given any $C'$ maximal $\dd$-independent subset of $A,$
$\card(C)=\card(C');$

(ii) any bijection $g: C\to C$ is elementary, moreover, if the kernel
is compact $g$ can be extended to an automorphism;

(iii) for finite $S\subs D$ \ $\card(\cl(S))\le \kappa;$

(iv) if $\lambda>\kappa$ then any definable subset of $D$ is either of
power  $\lambda$ or at most  $\kappa.$ \et
\pf (i) By Lemma~\ref{L1} $C'$ is in $\cl(C),$ the $\dd$-closure of $C.$ By
Lemma~\ref{star} $\card(C')\le \card(C).$ On the other hand $\cl(C')$  must
contain $C$ for otherwise $C'$ is not maximal. Hence $\card(C')=\card(C).$   \\

(ii) since $g$ is a (partial) isomorphism we may identify elements of
$C$ with their images to apply Proposition~\ref{p2}. Indeed, we can
see that q-construction $\{ a_j: j<\lambda\}$ of $\D$ as a q-prime
structure over $C$ gives
rise to the q-construction $\{ a_j: j<\lambda\}$ over same $C$ with
 types $\tau_j$ replaced by $g(\tau_j)$  obtained by exchanging
parameters in correspondence with $g.$ So by Proposition~\ref{p2} $g$
can be extended to an automorphism.    

(iii) By the definitions and the above proved there is a finite subset
$C_0\subs C,$ such that $S\subs \cl(C_0).$ Consider $E(C_0),$ which may be
assumed embedded in $E(C).$ Now construct inside $E(C)$ structure
$A=E(E(C_0)\cup C).$ Notice, that since $E(C_0)\subs \cl_A(C_0),$ \
the set $C\smin C_0$ is independent over $E(C_0).$\\
We now claim $\cl_A(C_0)=E(C_0).$
Let $a\in \cl_A(C_0).$ Since $A$ is q-atomic over $E(C_0)\cup C,$ there
is $\bar a$ extending $a,$ $\delta(\bar a/E(C_0)\cup C)=0$ and such that  $\bar a$ is q-atomic over
$E(C_0)\cup C.$ Let $\tau$ be 
 the q-atom of $\bar a$ over the set. More precisely, 
$$\tau=(V(\bar c), W(\ex(\bar c)),\{ W^{\frac{1}{l}}(\ex({\frac{1}{l}}\bar c)):l\in \N\})$$
with  $\bar c$ a finite
string from $C\smin C_0$  and
$(V,W)$ defined over a finite $B\subs E(C_0).$
On the other hand let
$$\tau_0=(V_0, W_0,\{ W_0^{\frac{1}{l}}:l\in \N\})$$
be the ex-locus ov $\bar a$ over $E(C_0).$
Since $\bar c$ is $\dd$-independent over $B\bar a,$
$\mdim(\bar a/B\bar c)=\mdim(\bar a/B)$ and hence $\dim V(\bar c)=\dim
V_0,$ $\dim W(\ex(\bar c))=\dim W_0.$ But $V_0$ and $W_0,$ being defined over
algebraically closed fields, are absolutely irreducible, hence
$\tau=\tau_0$
is defined over $B.$

 Since $E(C_0)$ is
existentially-algebraically closed, $\tau$ is realised in $E(C_0)$ by a tuple
$\bar a'.$ Since $\tau$ is a q-atom over $E(C_0)\cup C,$ we get $\bar
a=\bar a',$ proving the claim.
It follows from the claim and Proposition\ref{L4} that  $\card(\cl_A(C_0)\le \card(E(C_0))\le \kappa.$ On the other hand
$E(C)$ is embeddable in $A$ and by this embedding $\cl(C_0)$ goes into
$\cl_A(C_0).$  It finishes the proof of (iii).\\
(iv)  Choose $C_0\subs C$
finite such that $S\subs \cl(C_0).$
If a set $X$ definable over finite $S\le \D$ is of power greater than
$\kappa,$
then by (iii) it contains an element $b\notin \cl(C_0).$ But $b\in
\cl(C),$ so we can find finite $C_1\sups C_0$ and $c\in C$ such that 
$b\in \cl(C_1c)\smin \cl(C_1).$ Hence $\cl(C_1b)=\cl(C_1c).$ Consider
bijections $g:C\to C$ which fix $C_1.$ It follows that the
automorphism of $\D$ induced   by such $g$ have the property that
$g_1(b)\neq g_2(b)$ whenever $g_1(c)\neq g_2(c).$ Thus there are
$\lambda$ elements of the form $g(b),$ which all are in $X.$ \qed

\ssn{Pseudo-analytic dimension in canonical structures}
Throughout this section $\D$ is a canonical structure over a
$\dd$-independent  set of power $\lambda>\aleph_0$ and the standard kernel.
We expand the language by naming all elements of the standard kernel. \\ \\
\df For a subset $S\subs D^n$ definable over some finite $C\le \D$
{\bf the pseudo-analytic dimension of} $S,$ denoted  $\adim S,$ is given as
$$\adim S=\max_{\bar a\in S}\dd(\bar a/C).$$
\bl \lb{an1}
(i) The definition of $\adim$ does not depend on the choice of parameters
$C\le D.$\\
(ii) $\adim S\ge 0$ for any non-empty $S.$\\
(iii) For non-empty $S$ \
$\adim S=0\mbox{ iff }\card S \le\kappa.$\\
(iv)  For non-empty $S$ \ $\adim S\ge m>0$ iff $S$ can be projected to
$D^m$ so that the complement to the project
ion in $D^m$ is of pseudo-analytic
 dimension less than $m$ or empty.\\
(v) For algebraic varieties $V\subs D^n,$ $W\subs R^n$ \
$$\adim V=\dim V,\ \ \ \adim \ln W=\dim W.$$ \el
\pf Immediate.\qed
\bp \lb{p3} Suppose $(V,W)$ is a normal free pair over $C\le \D.$
Then $$\adim (V\cap \ln W)\ge \dim V+\dim W - n.$$\ep
\pf By induction on $d=\dim V+\dim W -n\ge 0.$

If $\dim W=n,$ then $V\cap \ln W=V.$ Thus we may assume $\dim V>1.$
Choose $c_1,\dots,c_n\in D$ such that $\dd(\la c_1,\dots,c_n\ra/C)=n.$
By Lemma~\ref{402}   we get then algebraic subvariety $V'\subs V$
  defined over
$C\cup\{ c_1,\dots,c_n\}$ such that pair $(V',W)$  is normal, free and
$$ \dim V'+\dim W'-n=d-1.$$ By induction  there is
$\bar a\in V\cap \ln W,$ \ $\dd(\bar a/C\bar c)=d-1.$ But, again by
Lemma~\ref{402}, $\dd(\bar c/C\bar a)<n=\dd(\bar c/C),$ hence
$\dd(\bar a/C)>\dd(\bar a/C\bar c)=d-1.$ \qed
\bl \lb{an2}
Let $(V,W)$ be a normal pair over $C\le \D$ and suppose
$$V\cap \ln W\neq\emptyset.$$
Then $$\adim(V\cap \ln W)\ge \dim V+\dim W-n.$$\el
\pf We use induction on $n.$ The case $n=1$ is trivial.\\
Consider the general case. Choose $\la a_1,\dots,a_{n}\ra\in V\cap \ln W.$
By the induction hypothesis there is
$\la a'_1,\dots,a'_{n-1}\ra \in V(a_n)\cap \ln W(\ex(a_n))$ with
$$\dd( \la a'_1,\dots,a'_{n-1}\ra/Ca_n)\ge \dim V(a_n)+\dim W(\ex(a_n))-(n-1).$$
If $ \dim V(a_n)=\dim V$ or $\dim W(\ex(a_n))=\dim W$ then
we would have the desired estimate
$$\dd( \la a'_1,\dots,a'_{n-1}\ra/Ca_n)\ge \dim V+\dim W-n.$$
Hence we assume   $ \dim V(a_n)=\dim V-1$ and $\dim W(\ex(a_n))=\dim W-1$
which yields  that the estimates are true for generic fibers and thus
the projections $V_n$ of $V$ and $W_n$ of $W$ to the last coordinate are
dense in $D$ and $R$ correspondingly.  Choose $a'_n\in D$ $\dd$-independent
over $C.$

By the induction hypothesis there is
$\la a'_1,\dots,a'_{n-1}\ra \in V(a'_n)\cap \ln W(\ex(a'_n))$ with
$$\dd( \la a'_1,\dots,a'_{n-1}\ra/Ca'_n)\ge \dim V(a_n)+
\dim W(\ex(a_n))-(n-1)=\dim V+\dim W-n-1.$$
By our choice
$$\dd( \la a'_1,\dots,a'_{n}\ra/C)=\dd( \la a'_1,\dots,a'_{n-1}\ra/Ca'_n)+
1\ge \dim V+\dim W-n$$ and $\la a'_1,\dots,a'_{n}\ra\in V\cap \ln W.$  \qed

{\bf References}\\

[B] S.Buechler {\em Essential Stability Theory}, Springer, Berlin 1996\\

[F] L.Fuchs {\em Infinite Abelian Groups}, Academic Press, New York -
London, 1970\\

[H] E.Hrushovski  {\em Strongly Minimal Expansions of Algebraically Closed Fields},
Israel Journal of Mathematics 79 (1992), 129-151               \\

[L] S.Lang {\em Introduction to Transcendental Numbers},
Addison-Wesley, Reading, Massachusets, 1966\\

[Sh] S.Shelah {\em Classification Theory},  revised edition, North-Holland--
Amsterdam--Tokyo 1990\\

[Z]  B.Zilber {\em Generalized Analytic Sets}, Algebra i Logika,
Novosibirsk,  v.36, N 4 (1997), 361 - 380

\end{document}